\documentclass[twoside, 10pt]{article}
\usepackage{cite}
\usepackage{amsmath}
\usepackage{paralist}
\usepackage{amssymb,mathrsfs,psfrag,graphicx}
\usepackage{hyperref}
\usepackage{refcheck}
\usepackage{amsfonts}
\usepackage{indentfirst}
\usepackage{latexsym}
\newtheorem{theorem}{Theorem}[section]
\newtheorem{remark}{Remark}[section]
\newtheorem{definition}{Definition}[section]
\newtheorem{lemma}[theorem]{Lemma}
\newtheorem{pro}{Proposition}[section]

\newcommand{\mr}{\mathbb{R}}

\newcommand{\na}{\nabla }

\newcommand{\dis}{\displaystyle}

\newcommand{\bl}{\begin{lemma}}
	\newcommand{\el}{\end{lemma}}
\newcommand{\et}{\end{theorem}}

\newcommand{\la}{\label}

\newcommand{\Om}{\Omega}

\newcommand{\bn}{\begin{eqnarray}}
\newcommand{\en}{\end{eqnarray}}
\newcommand{\bnn}{\begin{eqnarray*}}
\newcommand{\enn}{\end{eqnarray*}}

\newcommand{\bnnn}{\begin{eqnarray*}}
\newcommand{\ennn}{\end{eqnarray*}}

\newcommand{\ba}{\begin{aligned}}
\newcommand{\ea}{\end{aligned}}
\newcommand{\be}{\begin{equation}}
\newcommand{\ee}{\end{equation}}

\def\norm[#1]#2{\|#2\|_{#1}}

\makeatletter
\@addtoreset{equation}{section}
\makeatother
\title{ On regularity  and stability results of global weak solutions to  two-dimensional Navier-Stokes equations with density-dependent viscosity and vacuum
\thanks{ Partially supported by National Nature Science Foundation of China (Grant Nos. 12361046, 12661045) and Guizhou Provincial Basic Research Program General Project (Grant No. QianKeHeJiChu MS[2026]456).}}
\author{Yanfang Peng$^a$,  Hongyu Wang$^b$
\thanks{
	Email addresses: pyfang2005@sina.com (Y.F. Peng);  {hyuwang\_a@163.com} (H. Y. Wang).}\\[3mm]
\small \emph{ a. School of Mathematical Sciences, Guizhou Normal University, Guiyang 550025, P. R. China} \\ \small \emph {b.  School of Mathematics and computer Sciences,
Nanchang University, Nanchang, 330031,  P. R. China}}

\date{ }

\begin{document}
\maketitle
\begin{abstract}
	For the initial and  initial-boundary value problems of   the two-dimensional (2D) incompressible Navier-Stokes equations   with  density-dependent viscosity and vacuum,  Lions' global weak solutions  are shown to possess additional regularity, via establishing some key a priori decay-in-time rates for  the approximate problem,      under suitable
    smallness condition on   the  density-dependent viscosity.
	In addition, large-time decay rates of the pressure and  gradient  of velocity are    obtained,
    which are even exponential in time for  bounded
domains or the periodic case.  Note that these  results hold without any smallness conditions   on the initial data, even   for initially vanishing density.
	\end{abstract}	

\section{Introduction}
The 2D  nonhomogeneous incompressible Navier-Stokes equations with density-dependent viscosity  read  as follows:
\begin{equation}
	\begin{cases}
		\displaystyle\rho_{t}+\operatorname{div}\left(\rho u\right)=0, \\
		(\rho u)_{t}+\operatorname{div}\left(\rho u \otimes u\right)-\operatorname{div}(2 \mu(\rho) d)+\nabla P=0,\\
		\operatorname{div} u=0,
		\label{11}
	\end{cases}
\end{equation}
where $t\geq0$ is time, $x\in\Omega\subset \mr^2$ is the spatial coordinates, and $\rho=\rho(x,t), u=(u^1(x,t),u^2(x,t))$, and $P = P(x,t)$ represent the density, velocity,  and  pressure, respectively.  The deformation tensor is given by
\[d=\frac{1}{2}[\nabla u+(\nabla u)^{T}],\]
and the viscosity $\mu = \mu(\rho)$  satisfies the following hypothesis:
\begin{equation}
	\mu \in C[0, \infty), \quad \mu(\rho) >0.
	\label{12}
\end{equation}
System \eqref{11} provides a standard model for the evolution of multiphase flows, consisting of several immiscible, incompressible fluids with different densities.
In addition,  system \eqref{11} is subject to the initial conditions
\begin{equation}
\rho(x, 0)=\rho_{0}(x),\quad \rho u(x, 0)=\rho_{0} u_{0}(x),	\label{13}
\end{equation}
and  three types of  boundary conditions:
\begin{description}
	\item (i) (Periodic problem) $\Om=\mathbb{T}^2$ is the two-dimensional torus in $\mr^2$ and
	\begin{equation}\label{15}
		(\rho, u) \text{ is periodic in } x \text{ of  period } \Om;
	\end{equation}
	\item(ii) (Dirichlet problem) $\Om$ is a connected bounded smooth domain in $\mr^2$ and
	\begin{equation}\label{14}	u=0\,\,\text{ on } \partial\Om;\end{equation}
	\item (iii) (Cauchy problem) $\Om=\mr^2$. The far-field setting is understood through the finite-mass and energy spaces
\begin{equation}\label{16}
\begin{aligned}
	&\rho\in C([0,T];L^1(\mathbb R^2)),\qquad
	\sqrt\rho u\in L^\infty(0,T;L^2(\mathbb R^2)),\\
	&\nabla u\in L^2((0,T)\times\mathbb R^2),\qquad T>0.
\end{aligned}
\end{equation}
\end{description}

There is a large amount  of literature devoted to the existence and large-time behavior of weak solutions to \eqref{11}. When the density is  a positive  constant (i.e.,  the homogeneous case), we refer to \cite{Leray_1934, Duchon_1999, Chemin_2011,  Lions_1996,Temam_2001} and the references therein.
In the nonhomogeneous case, we briefly recall some well-posedness results in multi-dimensional spaces. For \eqref{11}  with constant viscosity, on the one hand, in the absence of vacuum, Kazhikov \cite{Kazhikov_1974} (see also \cite{AK_1973}) established the global existence of weak solutions. Later,
Ladyzhenskaya-Solonnikov \cite{Ladyzhenskaya_1978} first studied the unique resolvability of  system \eqref{11}  with   Dirichlet boundary conditions.
For  almost-critical Sobolev spaces,
Danchin \cite{Dan_2004} derived the local well-posedness in $\mr^n$ or the torus  $\mathbb{T}^n$ for any  $n$,  and the global one if $n=2$. Moreover,
   Global solutions to  three-dimensional
     Cauchy problem are also established   with small initial data in certain  scale invariant spaces (see\cite{AG_2011, AG_2012, Paicu_2012, AP_2007}) .  Under Dirichlet boundary condition,  Danchin-Mucha \cite{Dan_2013} investigated  a unique local  solution   with a piecewise constant density in dimension $n=2,3$. While
    the density is bounded above and below by positive constants,  Paicu-Zhang-Zhang  \cite{pzz_2013} constructed the global  well-posedness of  \eqref{11} in $\mr^n$ for dimensions $n=2,3,$ requiring    $ u_0\in  H^s(\mr^2)$ for some $s > 0$ ($n=2$) or $u_0\in   H^1(\mr^3)$ with $\|u_0\|_{L^2}\|\na u_0\|_{L^2}$ sufficiently small ($n=3$).
   Recently, Chen-Zhang-Zhao \cite{czz_2016} extended the 3D result in \cite{pzz_2013} to $u_0\in H^s(\mr^3)$ with $s>1/2$. Later, Zhang \cite{zhang_2020} obtained  global existence (see also \cite{Dan_2023}  for the uniqueness result) of solutions to  the three-dimensional 
    system \eqref{11}, with
 $\rho_0$ of a positive upper and lower bound  and  $u_0$ of small critical regularity.
In the  recent work of \cite{AG_2023,  AG_2021}, the authors proved  the global existence of unique solutions in $\mr^2$ for large initial data in critical spaces.  For the initial density allowed to vanish and constant viscosity,
   Simon \cite{Simon_1990} showed the global existence of weak solutions with finite energy.
    In \cite{Huang_2013}, Huang-Wang   pointed out that the unique local strong  solution obtained by  \cite{Choe_2003} was indeed global in  2D bounded domains. Furthermore, L\"{u}-Shi-Zhong \cite{Lu_2018} successfully obtained    the global existence and large-time asymptotic behavior of   solutions   in $\mr^2$ space by using
     the duality between BMO and Hardy spaces $\mathcal{H}^1$.
 Later, Danchin-Mucha \cite{dan_2019}   observed that while  $\rho_0$ is merely bounded and $u_0\in H^1(\Om)$, the global existence and uniqueness of  solutions hold in  bounded domains or  torus with general data if   $n=2$  and the      small scaling-invariant velocity if $n=3.$  Their results are extended to the   case of unbounded domains by  Prange-Tan \cite{jin_2024} with  compatibility condition
and some technical   conditions. Very  recently, Hao-Shao-Wei-Zhang \cite{wei_2025}  confirmed
    the global  weak solutions to the 2D  system \eqref{11} with  $0\leq\rho_0\in L^\infty(\mr^2)$ and $u_0\in L^2(\mr^2)$.Furthermore, they obtained the   uniqueness  when  $\rho_0$ remains  bounded away from zero.

For the case where   viscosity $\mu(\rho)$ depends on the density
$\rho$,  system  \eqref{11} becomes more complicated due to the strong coupling between viscosity and density.
On the  one hand, if   the initial density is strictly away from vacuum,  Abidi-Zhang \cite{AZ1_2015}  proved the global well-posedness   provided that  $\|\mu(\rho_0)-1\|_{L^\infty}\ll 1$ in $\mr^2$, and later by  \cite{AZ2_2015} if  $\|u_0\|_{L^2}\|\nabla u_0\|_{L^2}\ll 1$ and $\|\mu(\rho_0)-1\|_{L^\infty}\ll 1$ in $\mr^3$, respectively. More recently,  under the  Dirichlet boundary condition,
Huang-Li-Zhang \cite{huang2024} demonstrated
  the global well-posedness of solutions to the three-dimensional inhomogeneous Navier-Stokes 
   equations, with  large initial data  and  viscosity $\mu(\rho)=\mu\rho^\alpha (\mu>0, \alpha>1)$.
On the other hand, as for the initial density allowed to vanish, 
   Lions \cite{Lions_1996} obtained  the global existence of  weak solutions     with large data. However,  the uniqueness and regularity of the weak ones  presented in \cite{Lions_1996} are still open even in   2D case.   Later, Desjardins  \cite{Desjardins_1997} established  the  enhanced regularity for   periodic weak solutions  with   suitable assumptions on the density-dependent viscosity:
\begin{pro}\la{p1}
	In the notation used below, let $\Om=\mathbb{T}^2$. Assume that $\mu(\rho)$ satisfies \eqref{12} and  the initial data $(\rho_0, u_0)$ satisfies
	\[0\leq\rho_0\in L^\infty(\mathbb{T}^2),\quad u_0 \in H^1\left(\mathbb{T}^2\right),\quad\operatorname{div}u_0=0.\]
	Then there exists a positive constant $\varepsilon$ independent of the initial data such that if
	\begin{equation}\la{11111}
		\inf _{c>0}\left\|\frac{\mu\left(\rho_0\right)}{c}-1\right\|_{L^{\infty}\left(\mathbb{T}^2\right)}<\varepsilon,
	\end{equation}
	the following regularity results hold for  some positive constant $p^*$ and any $T>0$,
	\begin{equation*}
		\begin{cases}
			u \in L^{\infty}\left(0, T ; H^1\left(\mathbb{T}^2\right)\right), \sqrt{\rho} u_t \in L^{2}\left((0, T) \times \mathbb{T}^2\right),\\
			\left(\rho, \mu(\rho)\right) \in L^{\infty}\left((0, T) \times \mathbb{T}^2\right) \cap C\left([0, T] ; L^p\left(\mathbb{T}^2\right)\right), \text{ for all } p\in[1,\infty),\\
			\left(\nabla\left(P+R_i R_j\left(2 \mu(\rho) d_{i j}\right)\right), \nabla(\hat{P} \otimes Q(2 \mu(\rho) d))_{i j}\right) \in L^{2}\left((0, T)\times \mathbb{T}^2\right),\\
			\nabla u \in L^2\left((0, T) ; L^p\left(\mathbb{T}^2\right)\right), \text{ for all } p\in[4,p^*),
		\end{cases}
	\end{equation*}	with the Riesz transforms and projections specified in  \eqref{lj01} below.
\end{pro}

Later,   the  global
  well-posedness  of  \eqref{11} in two-dimensional bounded domains
     was established by Huang-Wang \cite{Huang_2014}   under  compatibility assumptions, provided that $\|\nabla \mu(\rho_0)\|_{L^q} (q>2)$ is small enough.   Moreover, under the small condition on $\|\nabla u_0\|_{L^2}$   and compatibility assumptions,  Huang-Wang \cite{Huang_2015}  demonstrated  the global existence of   solutions in three-dimensional
     bounded domains. Notably, in a breakthrough work He-Li-L\"{u} \cite{Li_2021}, the authors  first obtained  the global existence and exponential stability of   solutions to the Cauchy problem  in $\mr^3$, for suitably
small initial velocity in some homogeneous Sobolev space. More recently,  by exploiting the 2D cancellation structure and density moment bounds, Fan-L\"{u}-Wang \cite{FLW_2026} showed exponential decay of the velocity in $\mathbb{R}^2$  under small viscosity variation and a finite second moment of the initial density.

The primary objective of this paper is, under the viscosity assumptions of \cite{Desjardins_1997} and with merely
$L^2$ initial velocity, to construct global weak solutions to \eqref{11} with additional positive-time regularity. Moreover, our results hold not only for periodic problems, but also for Dirichlet and Cauchy problems.

Prior to presenting the   main results, we first explain some notation and conventions used
throughout this paper and clarify the precise meaning of weak solution.  We denote  $C$ to represent a positive generic constant which may vary in different cases, and  $C(p)$ to emphasize that $C$ depends on $p$. Moreover,  we set	$\varsigma(t)\triangleq\min\{1,t\},  B_N \triangleq\left\{x \in \mathbb{R}^2\mid | x|<N\right\}$ for $N>0$,  and  $R_i=\partial_i(-\Delta)^{-1/2}, (i=1,2),$  on  $\mathbb{T}^2$ or  $\mathbb{R}^2$,
with Fourier multiplier $i\xi_i/|\xi|$ for $\xi\ne0$. We use $(\nabla u)_{ij}=\partial_j u^i, (i,j=1,2) $    and \begin{equation}\label{lj01}
	\hat P=I-Q,\quad
	(Qf)_i=\partial_i\Delta^{-1}\partial_jf_j=-R_iR_jf_j.
\end{equation}	
Here $\Delta^{-1}$ has multiplier $-1/|\xi|^2$ at nonzero frequencies. On the torus $Q=0$ and $\hat P=I$ on constant vectors, and $\Delta^{-1}$ takes zero-mean inputs to zero-mean solutions. Thus $\hat P$ is the orthogonal projection onto divergence-free fields, including constant vectors in the periodic case. On a bounded Dirichlet domain, $\hat P$ denotes the $L^2$ orthogonal projection onto the space $L^2_\sigma(\Omega)$ introduced in the displayed space definitions below. On $\mathbb T^2$ or $\mathbb R^2$, for a matrix $M$, the tensor notation means 
\begin{equation}\label{lj02}
(\hat P\otimes QM)_{ij}
=\hat P_{ik}Q_{j\ell}M_{k\ell}
=\hat P_{ik}\partial_j\Delta^{-1}\partial_\ell M_{k\ell},
\quad i,j,k,\ell\in\{1,2\}, 
\end{equation}where repeated indices are summed and each multiplier acts on the indicated scalar component.  
 For  $k\geq1$ and $1 \leq p\leq\infty$, we introduce the standard Lebesgue and Sobolev spaces  as follows:
\begin{equation*}
	\left\{\begin{array}{l}
		\|f\|_{L^p}=\|f\|_{L^p(\Omega)},\quad \|f\|_{W^{k,p}}=\|f\|_{W^{k,p}(\Omega)},\\[3mm]\|f\|_{H^k}=\|f\|_{W^{k,2}(\Omega)},\quad C_{0, \sigma}^{\infty}=\left\{f \in [C_{0}^{\infty}(\Om)]^2:  \operatorname{div} f=0 \right\},\\[3mm]
L^2_{\sigma}=\overline{C_{0, \sigma}^{\infty}(\Om)} \quad \text{closure in the norm of } L^{2},\\[3mm]H_{0}^{1}=\overline{C_{0}^{\infty}(\Om)}\quad\text{closure in the norm of } H^{1},\\[3mm]  H_{0, \sigma}^{1}=\overline{C_{0, \sigma}^{\infty}(\Om)} \quad \text{closure in the norm of } H^{1},\\[3mm]	\tilde{D}^{1,2}=\tilde{D}^{1,2}(\Om)=\left\{v \in L_{\mathrm{loc}}^1\left(\Om\right) \mid \nabla v \in L^2\left(\Om\right)\right\},\\[3mm]\tilde{D}_{0, \sigma}^{1, 2}=\tilde{D}_{0, \sigma}^{1, 2}(\Om)=\overline{\left\{v \in [C_0^{\infty}(\Omega)]^2 \mid \operatorname{div} v=0\right\}} \text { closure in the norm of } \tilde{D}^{1, 2}.
	\end{array}\right.
\end{equation*}

We  consider  \eqref{11} and \eqref{12}  with the initial data $(\rho_0, u_0)$ satisfying
\begin{equation}\label{17} 0\leq\rho_{0}\leq\bar{\rho} \text{ a.e. in } \Om,  \quad \rho_0 \not \equiv 0, \quad
u_0\in L^2_\sigma,
\end{equation}
  In particular, in the periodic problem,  we may extend  $\rho_{0}$ and $u_0$ periodically on $\mathbb R^2$, and  if  $\Om=\mr^2$, we additionally assume $\rho_0\in L^1$.
Moreover, we denote 
\begin{equation} \label{ljo20}
	\underline{\mu}\triangleq \min_{\rho\in [0,\bar{\rho}]}\mu(\rho),\quad  \bar{\mu}\triangleq \max_{\rho\in [0,\bar{\rho}]}\mu(\rho).
\end{equation}

\begin{definition}\la{d1}{\rm We say that $(\rho, u, P)$ is a global weak solution to system \eqref{11} with initial data \eqref{13} and boundary condition \eqref{15}, or \eqref{14}, or \eqref{16} ,  if for all $T>0$,
		\begin{description}
			\item(i)  the density $\rho\geq0$ satisfies
			\[\rho\in L^\infty(\Om\times(0,\infty))\cap C([0,\infty);L^p), \,\, \text{for any}\, 1\leq p<\infty,\]
and  the velocity $u$ satisfies
\[
\begin{cases}
u\in L^2(0,T;H_{\mathrm{per}}^1)^2,&\text { for  periodic problem \eqref{15}},\\
u\in L^2(0,T;H^1_{0,\sigma}(\Omega))^2,& \text { for  Dirichlet problem  \eqref{14}},\\
\nabla u\in L^2(\Omega\times(0,T)),\quad
u\in L^2(0,T;H^1_{\rm loc}(\mathbb R^2))^2,& \text { for  Cauchy problem \eqref{16}}.
\end{cases}
\] Here $H_{\mathrm{per}}^1=\left\{u\in H^1(B_R) \text{ for all }  R<\infty,  u \text{ is periodic in } x_i \text{ of period } \Om \right\}$. In all three cases,
\[
\rho|u|^2\in L^\infty((0,\infty);L^1),\quad
\hat P(\rho u)\in C([0,\infty);L^2_w);
\]
		\end{description}
		\begin{description}
			\item(ii) $(\rho, u)$ solves  $\eqref{11}$  in the sense of distribution, that is, for all test functions $\varphi\in\Phi$ with compact support in time,
			$$
			\int_{\Om}\rho_0\varphi(\cdot,0)dx+\int_0^\infty\int_{\Om}(\rho\varphi_t+\rho u\cdot\na\varphi) dxdt=0,
			$$
			and
			$$\int_0^\infty\int_{\Om}u\cdot\na\varphi dxdt=0,$$
			where
			\[
\begin{aligned}
				&\Phi \triangleq \begin{cases} C_{\mathrm{per}}^{\infty}(\mathbb{R}^2 \times[0, \infty)) & \text { for  periodic problem \eqref{15}},\\C_0^{\infty}(\Omega \times[0, \infty)) & \text { for  Dirichlet problem  \eqref{14} or Cauchy problem \eqref{16}} \end{cases}
		\end{aligned}
		\]
			with
			$$	C_{\mathrm{per}}^{\infty}(\mathbb{R}^2 \times[0,\infty))=\left\{\varphi \in C^{\infty}(\mathbb{R}^2 \times[0, \infty)) \mid \varphi \text { is periodic in } x \text { of period } \Om \right\}.
		$$
			Moreover, 	for any  test functions  $\psi \in \Psi$ with compact support in time,
			$$
			\int_{\Om} \rho_0u_0 \cdot \psi(\cdot, 0) d x+\int_0^\infty\int_{\Om}\left(\rho u\cdot\psi_t+\rho u\otimes u:\na\psi\right) dxdt= \int_0^\infty\int_{\Om}2\mu(\rho)d:\na\psi dxdt, $$
			where
			$\Psi \triangleq\left\{ \psi \in \Phi^2 \mid \operatorname{div} \psi=0\right\}$;
		\end{description}	
		
		\begin{description}
			\item(iii) $\left(\rho, u\right)$  satisfies the energy inequality
			\begin{eqnarray*}
				\left(\frac{1}{2}\int_{\Omega} \rho\left|u\right|^2 d x\right)(t)+\int_0^t \int_{\Omega} 2\mu(\rho) \left(d:d\right) d x d s \leq \frac{1}{2}\int_{\Omega} {\rho_0}|u_0|^2 dx,
			\end{eqnarray*}
			for a.e. $ t \in(0, \infty).$
	\end{description}	}
\end{definition}

Our main results can be stated as follows:
\begin{theorem}\label{tm1}(Periodic problem) Let $\Om=\mathbb{T}^2$.
	Suppose that  \eqref{ljo20} and the initial data $\left(\rho_{0}, u_0\right)$  satisfy \eqref{17}.
	Then
	there exists a positive constant $\varepsilon_1$ depending only on the spatial dimensions   such that if
	\begin{equation}\label{110}
		{\inf_{c>0}\left\|\frac{\mu(\rho_0)}{c}-1\right\|_{L^{\infty}}}<\varepsilon_1,
	\end{equation}	
	the periodic problem \eqref{11}-\eqref{15}  has a global weak solution $\left(\rho, u, P\right)$  (with $\int_\Om P dx=0$)  satisfying that for any $0<\tau<T<\infty$ and  $p\in[1,\infty)$,
	\begin{equation}
		\begin{cases}\label{111}
			u \in L^{\infty}\left(\tau, T ; H^1\right), \sqrt{\rho} u_t \in L^2\left((\tau,T)\times\mathbb{T}^2\right)\cap L^{\infty}\left(\tau, T ; L^2\right),\\
			\left(\rho, \mu(\rho)\right) \in L^{\infty}\left((0, T) \times \mathbb{T}^2\right) \cap C\left([0, T) ; L^p\right),\\
			\nabla\left({P}+R_i R_j\left(2 \mu(\rho) d_{i j}\right)\right) \in L^2\left((\tau,T)\times\mathbb{T}^2\right)\cap L^{\infty}\left(\tau, T ; L^2\right),\\
			 \nabla(\hat{P} \otimes Q(2 \mu(\rho) d))_{i j} \in L^2\left((\tau,T)\times\mathbb{T}^2\right)\cap L^{\infty}\left(\tau, T; L^2\right).
		\end{cases}
	\end{equation}
	Moreover, for all  $p \in\left[4, p^*\right)$ and  some generic positive constant $C_0$ only depending on the spatial dimensions,
	\begin{equation}\label{1100}
		P \text{ and }\nabla u \in L^2\left((\tau, T); L^p\right) \cap L^\infty\left((\tau, T); L^p\right),	
	\end{equation}
	where
	\begin{equation}\label{114}
		\frac{1}{p^*}=2 C_0 \inf _{c>0}\left\|\frac{\mu(\rho_0)}{c}-1\right\|_{L^\infty},
	\end{equation}
		with $p^*=\infty$ when this infimum is zero, and
	\begin{equation}	
		\begin{aligned}&
			\|\nabla u(\cdot, t)\|_{L^{p}}  +\| P(\cdot, t)\|_{L^p}   \leq C  e^{-\frac\alpha2 t}, \quad  \text{for}\, \,\, t\ge 1,
		\end{aligned}\la{11110}
	\end{equation}
	with positive constants $\alpha$ and $C$,    depending only on    $\bar\rho $, $\underline{\mu}$, $\mathbb{T}^2$,	and  on $p$, $\bar\rho,$  $\bar{\mu},$ $\underline{\mu}$, $\|\rho_0\|_{L^1}^{-1},$  $\left\|u_0\right\|_{L^2}$,     $\mathbb{T}^2$, respectively.
\end{theorem}	

\begin{theorem}\label{tm3}(Dirichlet problem)
	Let $\Om\subset\mr^2$ be a connected smooth
	bounded domain.  Assume that  \eqref{ljo20} and  the initial data $\left(\rho_0, u_0\right)$ satisfy \eqref{17}. Then there exists a positive constant $\varepsilon_2$ depending only on $\Om$ and the spatial dimensions   such that if
	\begin{equation}\label{115}
		{\inf_{c>0}\left\|\frac{\mu(\rho_0)}{c}-1\right\|_{L^{\infty}}}<\varepsilon_2,
	\end{equation}
	the Dirichlet problem \eqref{11}-\eqref{13}, \eqref{14} has a   global weak solution $\left(\rho, u, P\right)$
    (with $\int_\Om P dx=0$)  satisfying that
	for any $0<\tau<T<\infty$ and $p\in[1,\infty)$,
	\begin{equation}
		\begin{cases}\label{116}
			u \in L^{\infty}\left(\tau, T; H^1\right), \sqrt{\rho} u_t \in L^2\left((\tau,T)\times\Om\right)\cap  L^\infty\left((\tau, T); L^2\right),\\
			\left(\rho, \mu(\rho)\right) \in L^{\infty}\left((0, T) \times \Om\right) \cap C\left([0, T) ; L^p\right),\\
			\nabla  P-\operatorname{div}\left(2 \mu(\rho) d\right) \in L^2\left((\tau,T)\times\Om\right)\cap L^\infty\left((\tau, T); L^2\right).
		\end{cases}
	\end{equation}
	Moreover, for all $ p \in\left[4, p^*\right),$
	\begin{equation}\label{117}
		P \text{ and } \nabla u \in L^2\left((\tau, T); L^p\right) \cap  L^\infty\left((\tau, T); L^p\right),
	\end{equation}	
		where there exists a continuous strictly decreasing function $\tilde C:(0,\varepsilon_2)\to(4,\infty)$, with the extended value $\tilde C(0)=\infty$, such that\begin{equation}\label{118}
		p^*=\tilde{C}\left({{\inf\limits_{c>0}}\left\|\frac{\mu(\rho_0)}{c}-1\right\|_{L^{\infty}}}\right),
	\end{equation}
	and
\begin{equation}\la{1161}	
		\begin{aligned}&
			\|\nabla u(\cdot,t)\|_{L^{p}}  +\| P(\cdot,t)\|_{L^p}   \leq C  e^{-\frac{\alpha_1} {2} t}, \quad  \text{for} \,\, t\ge 1,
		\end{aligned}
	\end{equation}
	with positive constants   $\alpha_1$ and $C$,  depending only on  $\bar\rho $, $\underline{\mu}$  $\Omega$,  and  on $p$, $\bar\rho,$  $\bar{\mu},$ $\underline{\mu}$,  $\left\|u_0\right\|_{L^2}$,   $\Omega,$ respectively.
\end{theorem}

\begin{theorem}\label{tm2}(Cauchy problem) Let $\Om=\mr^2$. Assume that  \eqref{ljo20} and  the initial data $\left(\rho_0, u_0\right)$ satisfy \eqref{17}. Then there is a positive constant $\varepsilon_3$ depending only on the spatial dimensions  such that if
	\begin{equation}\label{1101}
		{\inf_{c>0}\left\|\frac{\mu(\rho_0)}{c}-1\right\|_{L^{\infty}}}<\varepsilon_3,
	\end{equation}	
	the Cauchy problem  \eqref{11}-\eqref{13}, \eqref{16} has a global weak solution $\left(\rho, u,  \na P\right)$   satisfying that for any $0<\tau<T<\infty$ and  $p\in[1,\infty)$,
	\begin{equation}
		\begin{cases}\label{1111}
			\na u \in L^{\infty}\left(\tau, T; L^2\right), \quad\sqrt{\rho}\dot{u} \in L^2\left((\tau,T)\times\mathbb{R}^2\right)\cap L^{\infty}\left(\tau, T; L^2\right),
			\\
			\rho \in L^{\infty}\left((0, T) \times \mr^2\right) \cap C([0, T) ; L^p),\quad
 \mu(\rho) \in L^{\infty}\left((0, T) \times \mr^2\right) \cap C([0, T) ; L^p_{\text{loc}}),\\
			\left(\nabla\left(P+R_i R_j\left(2 \mu(\rho) d_{i j}\right)\right), \nabla(\hat{P} \otimes Q(2 \mu(\rho) d))_{i j}\right) \in L^2\left((\tau,T)\times\mathbb{R}^2\right)\cap L^{\infty}\left(\tau, T; L^2\right).
		\end{cases}
	\end{equation}
	Moreover,  for all $ p \in\left[4, p^*\right)$,
	\begin{equation}
	 	\nabla u \in L^2\left((\tau, T); L^p\right) \cap  L^\infty\left((\tau, T); L^p\right),
	\end{equation}	
	where there exists a continuous strictly decreasing function $\bar C:(0,\varepsilon_3)\to(4,\infty)$, with the extended value $\bar C(0)=\infty$, such that
	\begin{equation}\label{1141}
		p^*=\bar{C}\left({\inf\limits_{c>0}\left\|\frac{\mu(\rho_0)}{c}-1\right\|_{L^{\infty}}}\right),
	\end{equation}
	and
	\begin{equation}	
		\label{11611}
		\begin{aligned}&
			\|\nabla u(\cdot,t)\|_{L^{p}}^2     \leq C  t^{-1}, \quad  \text{for}\, any \,\, t\ge 1,
		\end{aligned}
	\end{equation}
	with the  positive constant   $C$  depending only on   $p$, $\|\rho_0\|_{L^1\cap L^\infty}$,  $\bar{\mu},$ $\underline{\mu}$,  and $\left\|u_0\right\|_{L^2}$.
\end{theorem}	

\begin{remark}
	Compared with some studies that  address the problem under  specific boundary condition, even for the constant viscosity case (see  \cite{dan_2019, AG_2021, AG_2023},) and the non-vacuum  case (see \cite{AZ1_2015, A_2007}), we solve the system \eqref{11}  under three different boundary conditions by a  decomposition approach   in Lemma \ref{lm22}. These results establish additional regularity  for the Lions'  weak solutions constructed here \cite{Lions_1996}.
\end{remark}

\begin{remark}
Eliminating  the initial data $\|\na u_0\|_{L^2}$ in Proposition \ref{p1}, Theorem \ref{tm1}    establishes the estimate on $L^\infty(\tau,T; L^p)$-norm of the gradient of the velocity and pressure, and  their  exponential decay-in-time rate for $0\leq\rho_0\in L^\infty$ and $u_0\in L^2$, which  advances the  global regularity  in \cite{Desjardins_1997}.
\end{remark}
\begin{remark}
As far as we can determine, 
 Theorem \ref{tm3}  gives a regularity result for  the Dirichlet problem to Navier-Stokes equations \eqref{11} with density-dependent viscosity and   vacuum. The key techniques are based on the classical Sobolev extension theorem and  the extension of BMO-Hardy space theory to bounded  domains.
\end{remark}
\begin{remark} Indeed, Theorem \ref{tm2} holds for
	the viscosity $\mu(\rho)$  satisfying  \eqref{1101} and arbitrarily large initial density (even
with compact support) with vacuum, which is in sharp contrast to \cite{AZ1_2015} and \cite{AZ2_2015} where they require smallness of $\|\mu(\rho_0)-1\|_{L^\infty}$ and the initial density to be strictly away from vacuum. \end{remark}

 \begin{remark}
 Nowadays, for constant viscosity, much research has further addressed the uniqueness of weak solutions to \eqref{11} in various domains, such as Danchin--Mucha \cite[Theorem 2.1]{dan_2019} in 2D bounded domains, Prange--Tan \cite[Theorem A]{jin_2024} for the Cauchy problem, Gancedo--Garc\'ia-Ju\'arez--Luna-Velasco
\cite[Theorem 1.4]{GancedoEtAl2025} for the periodic density-patch case. However, for variable viscosity, uniqueness remains open, and we will pursue this in future work.

\end{remark}

We now make some comments on the analysis of this paper.
To establish the  regularity  of global weak solutions of the initial and initial value boundary problems, it is vital to obtain time-independent  a priori estimates on smooth solutions
 for the approximated system (see \eqref{44}). As demonstrated in \cite{Desjardins_1997, Lu_2018, Huang_2014}, the key point lies in obtaining the bound on the $L^\infty(0, T;L^2)$-norm of $t^\frac{1}{2}\na u$.  In \eqref{44}, the estimates must be uniform in $\eta$ for the independently transported coefficient $h$.
To this end, we mainly employ some ideas from  \cite{Hoff_1995, Huang_2012, Li_2019, Lu_2018}. First,
we multiply $\eqref{44}_{3}$ by  $\dot{u}$ (with $\dot{u}\triangleq u_t+v\cdot\na u$) instead of the usual $u_t$ (see \cite{Desjardins_1997, Huang_2014}). It turns out that the crucial point is to bound the term (see \eqref{3100})
$$
G \triangleq \int_{\Om} P\partial_{i}v^{j}\partial_{j}u^{i}dx.
$$
Based on the following observation:
$$-c\Delta u+\na P =\operatorname{div}\left(2\left(h-c\right)d\right)-\rho\dot{u},$$
we  introduce a key Stokes  decomposition on bounded  or unbounded domains to derive integrability properties for $\na u$ and $P$ (more details are provided in
   \eqref{EP1} and \eqref{EP2}):
$$-c\Delta u_1+\na P_1=-\rho\dot{u},$$
and
$$-c\Delta u_2+\na P_2 =\operatorname{div}\left(2\left(h-c\right)d\right).$$
With this decomposition in hand,  $G$ is divided into  the following two parts:
$$G_1\triangleq\int_{\Om} P_{1}\partial_{i}v^{j}\partial_{j}u^{i}dx,
$$
and
$$
G_2\triangleq\int_{\Om} P_{2}\partial_{i}v^{j}\partial_{j}u^{i}dx.
$$

For $G_1$, in the periodic case, motivated by \cite{Lu_2018}, we  utilize  the duality of BMO and $\mathcal{H}^1$ to bound $G_1$ by   $\|\nabla u\|_{L^2}$ and $\|\sqrt{\rho}\dot{u}\|_{L^2}$.  However, the Dirichlet boundary case is more complicated  and requires extensions before applying the duality of BMO and $\mathcal{H}^1$. So, extending $u,v$ by zero to $\tilde{u}, \tilde v$, and $P_1$ to $\tilde{P_1}$ by the classical Sobolev extension  (see Extension Theorem in \cite{Evans_2010}),  we succeed in bounding $G_1$ (see Lemma \ref{lemmaP1}):
\[
\begin{aligned}
|G_1|&\leq C\|\tilde P_1\|_{\operatorname{BMO}(\mathbb R^2)}
 \|\partial_i\tilde v^j\partial_j\tilde u^i\|_{\mathcal H^1(\mathbb R^2)}\\
&\leq C\|P_1\|_{H^1(\Om)}\|\nabla v\|_2\|\nabla u\|_2
 \leq C\|P_1\|_{H^1(\Om)}\|\nabla u\|_2^2,
\end{aligned}
\]
where the last bound follows from \eqref{43}.

For $G_2$, the estimates in the periodic and Dirichlet cases are written in
\eqref{310} and \eqref{350}; both give the explicit bound
\[
 |G_2|\le \frac14\|\sqrt{\rho}\dot u\|_{L^2}^2
       +C\|\nabla u\|_{L^2}^4 .
\]
Next, to establish further regularity information, similarly to
\cite{Hoff_1995, Huang_2012, Lu_2018}, we apply the operator
$\partial_t+v\cdot\nabla$ to the $i$th component of $\eqref{44}_3$ and multiply by $\dot u^i$. This gives the $L^\infty_tL^2_x$ bound for $t\sqrt\rho\dot u$ and the $L^2_{t,x}$ bound for $t\nabla\dot u$ at small times (see \eqref{326}, \eqref{352}, \eqref{340}), and hence the $L^\infty_tL^p_x$ bounds for $t\nabla u$ and $tP$. The global weights are $\varsigma$ with an exponential factor in the periodic and Dirichlet cases, while the Cauchy estimates use the actual powers of $t$. In particular, the time weights allow merely $u_0\in L^2$.
Finally, with these a priori estimates at hand, we can establish global weak
solutions with the stated regularity. We also obtain the large-time decay
rates of the pressure and the velocity gradient, together with the exponential
rates in the Dirichlet and periodic cases (see \eqref{36} and \eqref{345}).
Moreover, it is necessary to  emphasize that for the whole space case, we  establish a priori estimates independent of $N$, by carrying out the same operation as the Dirichlet problem with $\Om=B_N$, and then obtain the expected   existence  and regularity of global weak solutions  by letting $N\rightarrow \infty$.

The rest of the paper is organized as follows. In Section \ref{2}, we list some known facts and elementary inequalities that will be used in  later analysis.  Section \ref{3} is devoted to establishing some necessary a priori estimates under three types of boundary conditions. Finally, we   construct approximate solutions and  pass the limit to complete  the proof of Theorems \ref{tm1}-\ref{tm2} in Section \ref{4}.

\section{Preliminaries}\label{2}
In this section, we present  some known facts and basic inequalities   that will be referenced later.

We start with the  following well-known Gagliardo-Nirenberg inequality  which will be used frequently later.
\begin{lemma}[Gagliardo-Nirenberg] \la{lm21} Let $\Om\subset\mr^2$ be a connected smooth
	bounded domain. Then   there exists a positive generic constant $C$,      depending on $\Om$ (specially  uniform for  $N$, if $\Om:=B_N$), such that for any $f \in H^{1}(\Omega)$ and  $p \in[2, \infty)$,
	\begin{equation}
		\begin{aligned}\label{21}
			\left\|f-\bar{f}\right\|_{L^{p}}\leq C\sqrt{p}\left\|f-\bar{f}\right\|_{L^{2}}^{2/p}\left\|\nabla f\right\|_{L^{2}}^{1-2/p},
		\end{aligned}
	\end{equation}
	where
	\begin{equation*}
		\bar{f} \triangleq  \frac{1}{\left|\Om\right|}\int_{\Omega}fd x.
	\end{equation*}
\end{lemma}
\noindent{\it Proof.} For the proof of \eqref{21}, we can refer to   \cite[Section 2]{Nirenberg_1959}. It remains to prove  the  uniformity for $N$, if  $\Om=B_N$.

Denote the rescaled function
$$f_N(x) \triangleq f(Nx), \quad x\in B_1.$$
Applying \eqref{21}  to $f_N$, we  see that
\begin{equation}\la{22}
	\left\|f_N-\bar{f}_N\right\|_{L^{p}\left(B_1\right)}\leq C\sqrt{p}\left\|f_N-\bar{f}_N\right\|_{L^{2}\left(B_1\right)}^{2/p}\left\|\nabla f_N\right\|_{L^{2}\left(B_1\right)}^{1-2/p}.
\end{equation}
Then by direct calculations
\begin{equation*}\begin{aligned}
		\int_{B_1}\left|f_N-\bar{f}_N\right|^pdx=&\int_{B_1}\left|f(Nx)-\frac{1}{\left|B_1\right|}\int_{B_1}f(Ny)dy\right|^pdx\\=&N^{-2}\int_{B_N}\left|f(x)-\frac{1}{\left|B_N\right|}
		\int_{B_N}f(y)dy\right|^pdx,
	\end{aligned}
\end{equation*}
and
\begin{equation*}
	\int_{B_1}\left|\na f_N\right|^2dx=N^2\int_{B_1}\left|\na f(Nx)\right|^2dx=\int_{B_N}\left|\na f\right|^2dx.
\end{equation*}
Substituting the two scaling identities into  \eqref{22}, we obtain  that there exists a positive constant $C$,  independent of $N$ such that
\begin{equation*}
	\left\|f-\bar{f}\right\|_{L^{p}\left(B_N\right)}\leq C\sqrt{p}\left\|f-\bar{f}\right\|_{L^{2}\left(B_N\right)}^{2/p}\left\|\nabla f\right\|_{L^{2}\left(B_N\right)}^{1-2/p},
\end{equation*}which completes the proof of Lemma \ref{lm21}. \hfill $\Box$

Let $\mathcal{H}^{1}(\mathbb{R}^2)$, $\mathcal{H}^{1}(\mathbb{T}^2)$ and BMO$(\mathbb{R}^{2})$, BMO$(\mathbb{T}^{2})$  denote the standard Hardy and $\mathrm {BMO}$ spaces, respectively.  In the periodic duality below, BMO is understood modulo constants and  paired with the zero-mean Hardy subspace.

Then, we state the following well-known facts which  are crucial to the proof in the next section (see    \cite[Theorem II.1(2) and Lemma II.1]{Clms},  \cite[Chapter IV, Section 1.2, Theorem 1]{Stein_1993}, \cite[Lemma 2.5]{Faraco_2021}).
\begin{lemma}\label{lemma2.3} Let $\Om=\mr^2$ or $\Om=\mathbb{T}^2$.

	{\rm(i)} There is a positive constant $C$ such that
	\begin{equation*}
		\|E \cdot B\|_{\mathcal{H}^{1}\left(\Om\right)} \leq C\|E\|_{L^{2}\left(\Om\right)}\|B\|_{L^{2}\left(\Om\right)},
	\end{equation*}	
	for all $E, B \in L^{2}\left(\Om\right)$ satisfying
	\begin{equation*}
		\operatorname{div} E=0,\quad \nabla^{\perp} \cdot B=0 \text { in } \mathcal{D}^{\prime}\left(\Om\right) .
	\end{equation*}
    Moreover, in the periodic case, $E$ and $B$ are also assumed to have zero mean.

	{\rm(ii)} There is a positive constant $C$ such that for any $v \in
	\tilde{D}^{1,2}(\Om)$,
	\begin{equation*}
		\|v\|_{\operatorname{BMO}\left(\Om\right)}  \leq C\|\nabla v\|_{L^{2}\left(\Om\right)},
	\end{equation*}
	 where the BMO semi-norm is taken modulo constants in the periodic case.
\end{lemma}

The following  commutator estimate  is used in the regularization argument of DiPerna–Lions    \cite{DiPernaLions_1989, Lions_1996}.
\begin{lemma}\la{lm26}
Let $\Omega\subset\mr^2$ be a domain and $f\in L^p(\Omega)$, $v\in W^{1,q}(\Omega;\mr^2)$  with $1<p,q<\infty$ and $1/p+1/q\le1$. For each compact $K\subset\Omega$ and sufficiently small $\varepsilon>0$,
\[
\|[\partial_i(fv)]_\varepsilon-\partial_i([f]_\varepsilon v)\|_{L^r(K)}
\le C\|f\|_{L^p(\Omega)}\|v\|_{W^{1,q}(\Omega)},
\]
where $C$ depends on $K,\Omega$ and the fixed mollifier, and $1/r=1/p+1/q$. If $\Omega=B_N$ and $v\in W^{1,q}_0(B_N;\mr^2)$, extend $f,v$ by zero to $\mr^2$, then
\begin{equation*}
\|[\partial_i(fv)]_\varepsilon-\partial_i([f]_\varepsilon v)\|_{L^r(B_N)}
\le C\|f\|_{L^p(B_N)}\| v\|_{W^{1,q}(B_N)},
\end{equation*}
with $C$ independent of $N$ and $\varepsilon$.
\end{lemma}

	Following the   method of \cite[Section 2.4]{Lions_1996}, we consider
    the following approximate problem
\begin{equation}\label{44}
	\begin{cases}
		\rho_t+\operatorname{div}\left(\rho v\right)=0, \\
		h_t+\operatorname{div}\left(hv\right)=0, \\
		(\rho u)_t+\operatorname{div}\left(\rho v \otimes u\right)-\operatorname{div}(2 h d)+\nabla P=0, \\
		\operatorname{div} u=0,
	\end{cases}\\	
\end{equation}
where $h=h(x,t)$ denotes the transported viscosity and the identity $h=\mu(\rho)$ is recovered only in the final constitutive limit.  Here $v$ is the common divergence-free transport velocity, whose precise definition is given in later sections. The problem \eqref{44} is subject to the  smooth initial conditions
\begin{equation}\label{46}
	\rho(x, 0)=\rho_0(x), \quad  u(x, 0)= u_0(x), \quad h(x,0)=h_0(x),
\end{equation} and the  boundary conditions below
\begin{description}
	\item (i) Periodic  case
	\begin{equation}\label{45}
		(\rho, u,h) \text{ is periodic in } x \text{ with  period } \Om;
	\end{equation}
	\item(ii) Dirichlet case for smooth bounded domain $\Om$
	\begin{equation}\label{450}
		u=0 \quad \text { on } \partial \Omega.
	\end{equation}
\end{description}
  Throughout the a priori estimates,  we assume  that $\underline\mu\le h_0\le\bar\mu$ and the small-oscillation hypothesis on $h_0$ reads as follows:
\begin{equation}\la{0}
    \inf_{c>0}\left\|\frac{h_0}{c}-1\right\|_{L^\infty}<\varepsilon.
\end{equation}
And we denote
\[
 \dot u\triangleq u_t+v\cdot\nabla u,\quad
 \dot v\triangleq\partial_tv+v\cdot\nabla v.
\]

And then, in view of the condition \eqref{0},   for fixed $\left(\rho, u, h,  P\right)$ (with $\int_{\Om}Pdx=0$) and any $c>0$, we introduce the following important decomposition for the approximated equations \eqref{44} as follows:

\begin{equation}
	\left\{
	\begin{aligned}
		-c\Delta u_{1}+\nabla P_{1}&=-\rho \dot{u},  &&x \in \Om, \\
		\operatorname{div} u_{1}&=0, && x \in \Om, \end{aligned}
	\right.
	\label{EP1}
\end{equation}
with $\int_\Om P_1dx=0$ and the periodic boundary condition
\begin{equation}\label{151}
	u_1 \text{ is periodic in } x \text{ of  period } \Om, \quad \int_\Om u_1dx=\int_\Om udx
\end{equation}
or  the Dirichlet boundary condition
\begin{equation}\label{141}
	u_1=0 \quad\text{ on } \partial\Om;\end{equation}
and
\begin{equation}
	\left\{
	\begin{aligned}
		-c\Delta u_{2}+\nabla P_{2}&=\operatorname{div}\left(2(h-c)d\right),&&  x \in \Om, \\
\operatorname{div}u_{2}&=0, &&  x \in \Om,
	\end{aligned}
	\right.
	\label{EP2}
\end{equation}
with $\int_\Om P_2dx=0$ and the periodic boundary condition
\begin{equation}\label{1511}
	u_2 \text{ is periodic in } x \text{ of  period } \Om, \quad \int_\Om u_2dx=0
\end{equation}
or  the Dirichlet boundary condition
\begin{equation}\label{1411}
	u_2=0 \quad\text{ on } \partial\Om.
\end{equation}

By the standard Stokes theory for linear systems, there exist unique solutions $(u_1, P_1)$ and $(u_2, P_2)$ to \eqref{EP1}, \eqref{151}, or \eqref{EP1}, \eqref{141}, and \eqref{EP2}, \eqref{1511}, or \eqref{EP2}, \eqref{1411}, respectively, such that
\begin{equation} \la{ljo8} u_1+u_2=u,  \quad P_1+P_2=P.
\end{equation}

Under the above decomposition \eqref{ljo8}, we have the following regularity estimates.
\begin{lemma}\label{lm22}
	Let  $\left(\rho, u, h, P\right)$ be a smooth solution to \eqref{44}-\eqref{45}, or \eqref{44}, \eqref{46}, \eqref{450}. Then for every $1<p<\infty$
and $c>0$, the  following estimates hold
	\begin{equation}
		\left\|c\nabla^{2} u_{1}\right\|_{L^p} +\left\|\nabla P_{1}\right\|_{L^p} \leq  {C}\left\|\rho \dot{u}\right\|_{L^p},
		\label{211}
	\end{equation}
	\begin{equation}
		\left\|c\nabla u_{2}\right\|_{L^{p}}+\left\|P_{2}\right\|_{L^{p}} \leq C  \left\| \left(h-c\right)d\right\|_{L^p}.
		\label{212}
	\end{equation}
    Here $C$ depends on $p$ and $\Om$, but not on $c$. For
$\Omega=B_N$, $C$ can be chosen independently of $N$.
\end{lemma}
\noindent{\it Proof.}
The standard $L^p$-estimate (see \cite[Section 1.2]{Temam_2001} and \cite[Section 4.1]{G_1998}) to  \eqref{EP1}, \eqref{151}, or \eqref{EP1}, \eqref{141}, and \eqref{EP2}, \eqref{1511}, or \eqref{EP2}, \eqref{1411}, yields \eqref{211} and \eqref{212} directly. We prove only the statement  about the uniformity for all domains $B_N$.

Indeed, let $y=N^{-1}x$, we have
\begin{equation*}
\begin{aligned}
-c\Delta_y[u_1(Ny)]+\nabla_y[NP_1(Ny)]
 &=-N^2(\rho\dot u)(Ny),
 \qquad y\in B_1.
\end{aligned}
\end{equation*}
where $\int_{B_1}P_1(Ny)dy=0$.
It follows from \eqref{211} that for any $p\in(1,\infty)$
\begin{equation*}
	\left\|c\nabla^{2} u_1(Ny)\right\|_{L^p(B_1)} +\left\|N\nabla P_{1}(Ny)\right\|_{L^ p(B_1)} \leq  {C}\left\|N^2(\rho\dot u)(Ny)\right\|_{L^p(B_1)},
\end{equation*}
together with transformation: $x=Ny$ gives \eqref{211} with  $C(p)$  independent of $N$. Similarly, \eqref{212} holds    uniformly  in   $N$. \hfill $\Box$

\section{A priori estimates}\label{3}
In this section, we will establish some necessary a priori estimates for  smooth  solution $(\rho, u, h,  P)$ to the approximated system \eqref{44}-\eqref{45}, or \eqref{44}, \eqref{46},  \eqref{450}. Let $T>0$ be a fixed time and $(\rho, u, h,  P)$ be the smooth solution to  \eqref{44}-\eqref{45}, or \eqref{44}, \eqref{46},  \eqref{450} on $\Om\times [0, T]$ with smooth initial data $(\rho_0, u_0)$ satisfying \eqref{17}.

\subsection{Periodic case}\label{3.1}
For $\eta>0$, we define  $ v\triangleq{u}*\omega_{\eta},$ where $\omega_{\eta}$ is the standard mollifying kernel of width $\eta$, and  also assume $\int_{\Om} \rho_0 u_0 d x=0$: this normalization is implemented at each approximation level and the original frame is restored in Section~\ref{4}.

Integrating the transport and momentum equations gives
\begin{equation}\label{periodic-conservation}
\int_\Omega\rho \,dx=\int_\Omega\rho_0\,dx>0,\quad
\int_\Omega\rho u\,dx=\int_\Omega\rho_0u_0\,dx=0,\quad
\int_\Omega\rho\dot u\,dx=0.
\end{equation}

\begin{lemma}\label{lm 30}  Let $(\rho,u,h,P)$ be a smooth solution to \eqref{44}-\eqref{45}. Then there exists some positive constant $C$ depending only on  $\bar\rho$,   $\underline{\mu}$,  $\left\|u_0\right\|_{L^{2}}$, $\Omega$ such that
	\begin{equation}
		\sup _{t \in[0, T]} e^{\alpha t} \left\|\sqrt\rho u\right\|_{L^2}^{2} +\int_{0}^{T} e^{\alpha t}\left\|\nabla u\right\|_{L^{2}}^{2} d t \leq C,
		\label{30}
	\end{equation}
	for some positive constant 	$\alpha$ depending only on  $\bar\rho $, $\underline{\mu}$, and $\Omega$.
\end{lemma}

\noindent{\it Proof.} According to standard arguments (see \cite[Theorem 2.1, Section 2]{Lions_1996}),
\begin{equation}
	0\le \rho \leq \bar \rho.
	\label{31}
\end{equation}

Multiplying $\eqref{44}_3$ by $u$ and integrating over $\Om$ leads to
\begin{equation}
	\frac{1}{2} \frac{d}{d t} \int_{\Omega} \rho\left|u\right|^{2} d x+2 \int_{\Omega} h\left|d\right|^{2} d x=0.
	\label{32}
\end{equation}
By the zero total momentum in \eqref{periodic-conservation},
\begin{equation*}
\int_\Omega\rho\left|u-\frac1{|\Omega|}\int_\Omega u\,dx\right|^2dx
=\int_\Omega\rho|u|^2dx+
\left(\int_\Omega\rho\,dx\right)\left|\frac1{|\Omega|}\int_\Omega u\,dx\right|^2,
\end{equation*}
which, combined with  \eqref{31} and   Poincar\'e inequality,  yields
\begin{equation*}
\begin{aligned}
\int_\Omega\rho|u|^2dx
&\le\int_\Omega\rho\left|u-\frac1{|\Omega|}\int_\Omega u\,dx\right|^2dx\\
&\le2\bar\rho\left\|u-\frac1{|\Omega|}\int_\Omega u\,dx\right\|_2^2
\le C(\Omega)\bar\rho\|\nabla u\|_2^2.
\end{aligned}
\end{equation*}
Using this, $\eqref{44}_2$, $\eqref{44}_4$,  \eqref{31},  $\underline\mu\le h_0\le\bar\mu$ and Poincar\'{e} inequality gives
\begin{equation}
	\int_{\Omega} \rho\left|u\right|^{2} dx \leq C\int_\Om\left|\na u\right|^2dx\leq\frac{1}{\alpha}\int_{\Om}h\left|d\right|^2dx,
	\label{33}
\end{equation}
for  some positive constant $\alpha$ depending only on $\bar{\rho}$ and $\Om$.  Multiplying this by $e^{\alpha t}$ and integrating  over $[0,T]$, one has
\be\nonumber
\sup _{t \in[0, T]} e^{\alpha t} \left\|\sqrt\rho u\right\|_{L^2}^{2} + \int_0^T e^{\alpha t}\left(\int_{\Omega}h\left|d\right|^{2}dx\right) dt\leq C,
\ee
which together  with the second inequality of  \eqref{33} gives \eqref{30}. For later use, if a periodic vector field $f$ satisfies $\int_\Omega\rho f\,dx=0$, then
\begin{equation*}
\left(\int_\Omega\rho\,dx\right)\frac1{|\Omega|}\int_\Omega f\,dx
=-\int_\Omega\rho\left(f-\frac1{|\Omega|}\int_\Omega f\,dx\right)dx.
\end{equation*}
 Poincar\'{e} inequality therefore implies
\begin{equation}\label{periodic-weighted-mean}
\|f\|_{L^2}\le C(\Omega)\left(1+\frac{\|\rho\|_{L^2}}{\int_\Omega\rho\,dx}\right)
\|\nabla f\|_{L^2},
\end{equation}
which together with the Gagliardo-Nirenberg inequality gives $\|f\|_{L^p}\le C\|\nabla f\|_{L^2}$ for any   $p\in[2,\infty)$, with $C$ depending on  $\Om$,  $\bar\rho$,  and  $\|\rho_0\|_{L^1}^{-1}$. Later, this result will be applied to $f=u$ and $f=\dot u$, respectively. The proof of Lemma \ref{lm 30} is completed.  \hfill $\Box$

From condition \eqref{0}, we now present some fundamental integrability results on $\na u$ and $ P$.
\begin{lemma}\label{lm24}
	Let  $\left(\rho, u,h, P\right)$ be a smooth solution to \eqref{44}-\eqref{45} satisfying \eqref{0}. Then  there exists  a positive constant  $\varepsilon$  depending only  on the spatial dimensions,   such that  for $p\in[4,p^*)$,
	\begin{equation}\label{210}
		\left\|\na u\right\|_{L^p}\le C\left\|\na u\right\|_{L^2}^{2/p}\left\|\rho \dot{u}\right\|_{L^2}^{1-2/p},
	\end{equation}
	and
	\begin{equation}\la{2101}
		\left\|P\right\|_{L^p}\leq C\left\|{\rho}\dot{ u}\right\|_{L^2}+C\left\|\na u\right\|_{L^2},
	\end{equation}
	where  $C$  depends only on $p$,  $\bar{\mu}$,  $\underline{\mu}$,  $\Om$, and for some universal constant $C_0>0$,
	$$	\frac{1}{p^*}=2 C_0 \inf _{c>0}\left\|\frac{h_0}{c}-1\right\|_{L^\infty}.$$
    Here we set $p^*=\infty$,  when the infimum is zero.
\end{lemma}
\noindent{\it Proof.}
By \eqref{211}, \eqref{212}, and \eqref{periodic-conservation}, the   Stokes equations \eqref{EP1} and \eqref{EP2} give  that for any $c>0$,
\begin{equation*}
c\partial_i\partial_ju_1=-R_iR_j\hat P(\rho\dot u),\qquad
\nabla u_2=-2(\hat P\otimes Q)\big((h/c-1)d\big).
\end{equation*}
Now,  we take
\begin{equation}\label{280}
\underline\mu\le c=c_0\triangleq \frac{\operatorname*{ess\,inf}h_0+\operatorname*{ess\,sup}h_0}{2}\le\bar\mu.
\end{equation}
Then, it follows from these identities   and the Gagliardo-Nirenberg inequality that for
any $p\in (2,\infty)$
\begin{equation}\label{28}
\begin{aligned}
\|\nabla u\|_{L^p}
&\le\|\nabla u_1\|_{L^p}+\|\nabla u_2\|_{L^p}\\
&\le C\sqrt p\,\|\nabla^2u_1\|_{L^2}^{1-2/p}\|\nabla u_1\|_{L^2}^{2/p}+\|\nabla u_2\|_{L^p}\\
&\le C\sqrt p\,c_0^{-(1-2/p)}\|\rho\dot u\|_{L^2}^{1-2/p}
\left(\|\nabla u\|_{L^2}+C\|(h/c_0-1)d\|_{L^2}\right)^{2/p}
+C_0p\|(h/c_0-1)d\|_{L^p}\\
&\le Cp^{3/2}c^{-(1-2/p)}(1+\|h/c_0-1\|_{L^\infty})
\|\nabla u\|_{L^2}^{2/p}\|\rho\dot u\|_{L^2}^{1-2/p}
+C_0p\|h/c_0-1\|_{L^\infty}\|\nabla u\|_{L^p}\\&\le Cp^{3/2}\underline{\mu}^{-(1-2/p)}(2+\bar{\mu}/\underline{\mu})\|\nabla u\|_{L^2}^{2/p}\|\rho\dot u\|_{L^2}^{1-2/p}
+C_0p\|h/c_0-1\|_{L^\infty}\|\nabla u\|_{L^p}\\
&= Cp^{3/2}\underline{\mu}^{-(1-2/p)}(2+\bar{\mu}/\underline{\mu})\|\nabla u\|_{L^2}^{2/p}\|\rho\dot u\|_{L^2}^{1-2/p}
+C_0p\inf_{c>0}\left\|\frac{h_0}{c}-1\right\|_{L^\infty}\|\nabla u\|_{L^p}
\end{aligned}
\end{equation}
where $C_0>0$  (only depending on the spatial dimensions) denotes the smallest constant $C>0$ such that  (see (56) in  \cite[Section 1.2]{Desjardins_1997})
$$
\forall p \in[2, \infty), \quad M\in C^{\infty}\left(\mathbb{T}^2 ; \mathcal{M}_{n \times n}(\mathbb{R})\right), \quad\left\|\hat P \otimes Q M\right\|_{L^p} \le \frac{C}{2} p\left\|M\right\|_{L^p},
$$
and in the last equality, we have used
$$\left\|\frac {h}{c_0}-1\right\|_{L^\infty}=\left\|\frac{h_0}{c_0}-1\right\|_{L^\infty}=\inf_{c>0}\left\|\frac {h_0}{c}-1\right\|_{L^\infty},$$
due to $\eqref{44}_2$ and the fact that  for any $c>0$,
\begin{align*}
\left\|\frac{h_0}{c}-1\right\|_{L^\infty}
&=\max\left\{
\left|\frac{\operatorname*{ess\,inf}h_0}{c}-1\right|,
\left|\frac{\operatorname*{ess\,sup}h_0}{c}-1\right|
\right\}\\
&=
\begin{cases}
\displaystyle\frac{\operatorname*{ess\,sup}h_0}{c}-1,
&0<c\le c_0,\\[6pt]
\displaystyle 1-\frac{\operatorname*{ess\,inf}h_0}{c},
&c>c_0
\end{cases}\\
&\ge
\frac{\operatorname*{ess\,sup}h_0-
      \operatorname*{ess\,inf}h_0}
     {\operatorname*{ess\,sup}h_0+
      \operatorname*{ess\,inf}h_0}
=\left\|\frac{h_0}{c_0}-1\right\|_{L^\infty}.
\end{align*}

Choosing $\varepsilon$ in \eqref{0} and $p*$  as follows:
$$\varepsilon=\frac1{8C_0}, \quad p*=\frac{1}{2C_0\inf\limits_{c>0}\left\|\frac{h_0}{c}-1\right\|_{L^{\infty}}},$$  one gets from   \eqref{280} and  \eqref{28} that  for $4\le p<p^*$,
\begin{equation}\label{281}
\|\nabla u\|_{L^p}
\le C\|\nabla u\|_{L^2}^{2/p}\|\rho\dot u\|_{L^2}^{1-2/p}
+\frac{p}{2p^*}\|\nabla u\|_{L^p}
\le C\|\nabla u\|_{L^2}^{2/p}\|\rho\dot u\|_{L^2}^{1-2/p}
+\frac12\|\nabla u\|_{L^p},
\end{equation}
which gives \eqref{210}. In particular, we set $p^*=\infty$, when  $\dis\inf_{c>0}\left\|\frac {h_0}{c}-1\right\|_{L^\infty}=0$.

Then, \eqref{210},  combined with \eqref{ljo8}-\eqref{212},  yields that,  for $p\in[4, p^*)$,
\begin{equation*}\begin{aligned}
		\left\|P\right\|_{L^p}\leq& \left\|P_1\right\|_{L^p}+\left\|P_2\right\|_{L^p}\\
		\leq&C\left\|\na P_1\right\|_{L^{2}}+C\left\|\na u \right\|_{L^p}\\\leq
		&C\left\|\rho\dot{ u}\right\|_{L^2}+C\left\|\rho\dot{u}\right\|_{L^2}^{1-2/p}\left\|\na u\right\|_{L^2}^{2/p}\\\leq
		&C\left\|{\rho}\dot{ u}\right\|_{L^2}+C\left\|\na u\right\|_{L^2}.
	\end{aligned}
\end{equation*}
Thus, we  complete the proof of Lemma \ref{lm24}.\hfill $\Box$

With Lemma \ref{lm24} at hand, we establish   the following exponential decay rate for $\|\nabla u\|_{L^2}.$
\begin{lemma}\label{lm31}
	Let $(\rho,u,h,P)$ be a smooth solution to \eqref{44}-\eqref{45} satisfying \eqref{0}. For $p\in[4,p^*)$,  there exists  a positive constant $C$  depending only on $\bar{\mu}$, $\underline{\mu}$, $\bar\rho$, $\|\rho_0\|_{L^1}^{-1}$,  $p$, $\left\|u_0\right\|_{L^2}$,  and $\Omega$ such that
	\begin{equation}
		\sup _{t \in[0, T]} \varsigma e^{\alpha t}\left\|\nabla u\right\|_{L^{2}}^{2}+\int_{0}^{T}  \varsigma e^{\alpha t}\left(\varsigma\left\|\sqrt{\rho}{u_t}\right\|_{L^2}^2+\left\|\sqrt{\rho}\dot{u}\right\|_{L^2}^2+\left\|\na u\right\|_{L^p}^2+\left\|P\right\|_{L^p}^2\right)dt \leq C,
		\label{36}
	\end{equation}
	where  $\alpha$ and  $\varepsilon$  are as  Lemmas \ref{lm 30} and  \ref{lm24}, respectively.
\end{lemma}
\noindent{\it Proof.}
By \cite[Appendix C.5, Theorem 7]{Evans_2010}, we have, for each fixed periodic $u\in W^{1,p}(\Om)$ and $1<p<\infty$,
\begin{equation*}
v\rightarrow u\quad\text{in }W^{1,p}(\Om)\quad\text{as }\eta\rightarrow0,
\end{equation*}
and
\begin{equation}\label{471}
\| v\|_{W^{1,p}(\Om)}\le\|u\|_{W^{1,p}(\Om)}.
\end{equation}

Integrating  $\eqref{44}_3$ multiplied  by $\dot{u}$  over $\Omega$ leads to
\begin{equation}
	\int_{\Omega} \rho|\dot{u}|^{2} dx=\int_{\Omega} \dot{u}\cdot\mathrm{div}\left(2hd\right)-\dot{u}\cdot \nabla P  dx \triangleq  I_{1}+I_{2} .
	\label{38}
\end{equation}
For $I_1$, it follows from    \eqref{210} (for $p=4$), \eqref{31} and \eqref{471}  that
\begin{equation}
	\begin{aligned}
		I_{1}=& \int_{\Omega}  \dot{u}\cdot\operatorname{div}\left(2 h d\right)  dx \\
		=&-\int_{\Omega} h \partial_{t}\left|d\right|^{2} dx-\int_{\Omega} h\left(\partial_{i} u^{j}+\partial_{j} u^{i}\right) \partial_{j}\left(v^{k} \partial_{k} u^{i}\right) dx \\
		=&-\frac{d}{dt}\int_{\Omega} h\left|d\right|^{2}dx+\int_{\Omega} \partial_th \left|d\right|^{2}dx-\int_{\Omega} hv^{k} \partial_{k} \left|d\right|^{2} dx \\
		& -\int_{\Omega} h\left(\partial_{i} u^{j}+\partial_{j} u^{i}\right) \partial_{j}v^k \partial_{k}u^idx \\
		=&-\frac{d}{dt}\int_{\Omega} h\left|d\right|^{2}dx-\int_{\Omega} h\left(\partial_{i} u^{j}+\partial_{j} u^{i}\right) \partial_{j}v^k \partial_{k}u^idx\\ \leq &-\frac{d}{dt}\int_{\Omega} h\left|d\right|^{2}dx+C\left\|\nabla v\right\|_{L^2}\left\|\na u\right\|_{L^4}^2\\
		\leq &-\frac{d}{d t} \int_{\Omega} h|d|^{2} dx +\frac14\|\sqrt{\rho} \dot{u}\|_{L^{2}}^{2}+C \left \| \nabla u\right \|_{L^{2}}^{4},
	\end{aligned}
	\label{39}
\end{equation}
where  in the fourth equality, we have used  \eqref{44}$_2$:
\begin{equation}\label{10000}
	\partial_th+\operatorname{div} \left(hv\right)=0.
\end{equation}
Moreover, utilizing  \eqref{ljo8} and periodic Hardy--BMO duality, we get
\begin{equation}
	\begin{aligned}		
		I_{2}
		&=-\int_{\Omega} \nabla P \cdot\left(\partial_{t} u+v \cdot \nabla u\right)dx\\
		&=\int_{\Omega} P \partial_{i} v^{j} \partial_{j} u^{i}dx\\
		& = \int_{\Omega} P_1 \partial_i v^j \partial_j u^i dx + \int_{\Omega} P_2 \partial_i v^j \partial_j u^i dx \\&\leq C\|{P}_{1}\|_{\text{BMO}}\left\|\partial_{i} {v}^{j} \partial_{j} {u}^{i}\right\|_{\mathcal{H}^{1}}+C\left\|P_{2}\right\|_{L^{4}}\left\|\nabla u\right\|_{L^{4}}\left\|\nabla v\right\|_{L^{2}}.
	\end{aligned}
	\label{3100}
\end{equation}
Using $\operatorname{div}\left(\partial_{ i}v\right)=\partial_{ i}\operatorname{div}v, \nabla^{\perp} \cdot\left(\nabla {u}^{i}\right) =0,$   Lemma \ref{lemma2.3} yields $$\|{P}_{1}\|_{\text{BMO}}\left\|\partial_{i} {v}^{j} \partial_{j} {u}^{i}\right\|_{\mathcal{H}^{1}}\leq C\|\na{P}_{1}\|_{L^2}\left\|\na v\right\|_{L^2}\left\|\na u\right\|_{L^2},$$
which, together with  \eqref{211}, \eqref{212}, \eqref{31}, \eqref{210}, \eqref{471}, and \eqref{3100}, shows that
\begin{equation}
	\begin{aligned}		
		I_{2}
		&\leq C\left\|\nabla v\right\|_{L^{2}}\left(\left\|\nabla P_{1}\right\|_{L^{2}}\left\|\nabla u\right\|_{L^{2}} + \left\|\nabla u\right\|_{L^{4}}^2\right).\\
		&\leq C\|\nabla u\|_{L^{2}}^{2}\|\sqrt{\rho} \dot{u}\|_{L^{2}}\\
		&\leq \frac14\|\sqrt{\rho} \dot{u}\|_{L^{2}}^{2}+C\|\nabla u\|_{L^{2}}^{4}.
	\end{aligned}
	\label{310}
\end{equation}

Now, inserting \eqref{39}  and \eqref{310} into \eqref{38} gives
\begin{equation*}
	\frac{d}{d t} \int_{\Omega} h \left|d\right|^{2} dx+\frac{1}{2}\int_{\Omega} \rho\left|\dot{u}\right|^{2} dx\le C\left\|\nabla u\right\|_{L^{2}}^4.
\end{equation*}
Multiplying this by $\varsigma e^{\alpha t}$, then applying \eqref{33} and Gronwall's inequality, we easily check that
\begin{equation}\la{322}
	\sup _{t \in[0, T]}  \varsigma e^{\alpha t}\left\|\nabla u\right\|_{L^{2}}^{2}+\int_{0}^{T}  \varsigma e^{\alpha t}\left\|\sqrt{\rho}\dot{u}\right\|_{L^2}^2dt \leq C,
\end{equation}
which, together  with  \eqref{31},  \eqref{210} and \eqref{2101}, leads to for $p\in[4,p^*)$\begin{equation}\label{323}
	\int_{\Om} \varsigma e^{\alpha t}\left(\left\|\nabla u\right\|_{L^{p}}^{2}+\left\| P\right\|_{L^{p}}^{2}\right) \leq C.
\end{equation}

Using \eqref{periodic-weighted-mean} and \eqref{471}, we have $\|v\|_4\le\|u\|_4\le C\|\nabla u\|_2$. Consequently,  \eqref{31} gives
\begin{equation}
	\begin{aligned}
		\int_{\Omega} \rho |u_{t}|^{2} d x & \leqslant \int_{\Omega} \rho|\dot{u}-v \cdot \nabla u|^{2} d x \\
		& \leq  C \left\|\sqrt{\rho} \dot{u}\right\|_{L^2}^2+C \int_{\Omega} \rho|v \cdot \nabla u|^{2} d x \\
		& \leq C \left\|\sqrt{\rho} \dot{u}\right\|_{L^2}^2+C\left \|v\right \|_{L^{4}}^2\left \|\nabla u\right \|_{L^{4}}^{2} \\
		& \leq C \left\|\sqrt{\rho} \dot{u}\right\|_{L^2}^2+C\left \|\nabla u\right \|_{L^{2}}^{2}\left\|\na u\right\|_{L^4}^2.
	\end{aligned}
	\label{3361}
\end{equation}
Hence, \eqref{322}-\eqref{3361} imply
\begin{equation}\label{periodic-ut-integral}
\begin{aligned}
&\int_0^T\varsigma^2e^{\alpha t}\|\sqrt\rho u_t\|_{L^2}^2dt\\
&\le C\int_0^T\varsigma e^{\alpha t}\|\sqrt\rho\dot u\|_{L^2}^2dt
+C\left(\sup_{0<t\le T}\varsigma\|\nabla u(t)\|_{L^2}^2\right)
\int_0^T\varsigma e^{\alpha t}\|\nabla u\|_{L^4}^2dt\le C.
\end{aligned}
\end{equation}
Combining  \eqref{322}-\eqref{periodic-ut-integral}  proves \eqref{36} and completes  the proof of Lemma \ref{lm31} .\hfill $\Box$

\begin{lemma}\label{lm32}
	Let $(\rho,u,h,P)$ be a smooth solution to \eqref{44}-\eqref{45} satisfying \eqref{0}. For any $p\in[4,p^*)$, 	  there exists a positive constant $C$ with the same dependencies as in Lemma \ref{lm31}, such that
	\begin{equation}
		\sup _{t \in[0, T]} \varsigma^2 e^{\alpha t}\left(\varsigma^{\frac{1}{2}}\left\|\sqrt{\rho} {u_t}\right\|_{L^{2}}^{2}+\left\|\sqrt{\rho} \dot{u}\right\|_{L^{2}}^{2}+\left\|\na u\right\|_{L^p}^2+\left\|P\right\|_{L^p}^2\right)+\int_{0}^{T} \varsigma^2 e^{\alpha t} \left\|\nabla \dot{u}\right\|_{L^{2}}^{2} d t \leq C,
		\label{326}
	\end{equation}
	where $\alpha$ and $p^*$ are   as Lemmas \ref{lm 30} and \ref{lm24} respectively.
\end{lemma}
\noindent{\it Proof.}
Operating $\left(\partial_t +v\cdot \nabla \right)$ to $\eqref{44}^{i}_3$ indicates  that
\begin{equation}
	\begin{aligned}
		\partial_{t}&\left(\rho \dot{u}^{i}\right)+\operatorname{div}\left(\rho v \dot{u}^{i}\right)-\operatorname{div}\left(2hD^{i}\right)\\
		&=-\partial_j\left(h\partial_j v^k \partial_k u^i\right)-\partial_j\left(h\partial_i v^k \partial_k u^j\right)\\
		&\;\;\;\;-\partial_j v^k \partial_k \left(h(\partial_i u^j+\partial_j u^i)\right)-\left(\partial_t +v\cdot \nabla \right)\partial_{ i}P
		,
	\end{aligned}
	\label{327}
\end{equation}
where
$$
D\triangleq\frac{1}{2}\left[\nabla \dot{u}
+\left(\nabla \dot{u}\right)^{T}\right]
$$ satisfies
\begin{equation}
	\begin{aligned}
		2\int_{\Omega}\left|D\right|^{2}\mathrm{~d}x
		&=\frac{1}{2}\int_{\Omega}\left(\partial_{i}\dot {u}^{j}+\partial_{j}\dot{u}^{i} \right)^{2}dx\\
		&=\int_{\Omega}\left|\nabla \dot{u}\right|^{2}dx+\int_{\Omega}\partial_{i}\dot{u}^{j}\partial_{j}\dot{u}^{i}dx\\
		&=\int_{\Omega}\left|\nabla \dot{u}\right|^{2}dx+\int_{\Omega}\left|\operatorname{div} \dot u\right|^2 dx.
	\end{aligned}
	\label{332}
\end{equation}

Then, integrating \eqref{327} multiplied by $\dot{u}^i$ over $\Omega$ leads to
\begin{equation}
	\begin{aligned}
		\frac{1}{2} &\frac{d}{dt} \int_\Omega \rho\left|\dot{u}\right|^{2} dx+ 2 \int_{\Omega} h\left|D\right|^{2} dx\\
		&=-\int_\Omega\dot{u}^i\partial_j\left(h\partial_j v^k \partial_k u^i\right)  dx-\int_\Omega \dot{u}^i\partial_j\left(h\partial_i v^k \partial_k u^j\right)  dx\\
		& \;\;\;\;-\int_\Omega \dot{u}^i\partial_j v^k \partial_k \left(h(\partial_i u^j+\partial_j u^i)\right) dx-\int_\Om\dot{u}^i \left(\partial_t +v\cdot \nabla \right)\partial_{ i}Pdx \\
		&\triangleq \sum_{i=1}^4 J_i.
	\end{aligned}
	\label{328}
\end{equation}

Each term on the right-hand side of \eqref{328} is estimated as follows.

First,  by \eqref{210} (for $p=4$),  \eqref{471},  and periodic integration by parts, one gets
\begin{equation}
	\begin{aligned}
		J_1+J_2+J_3= &\int_\Omega h\partial_j v^k \partial_k u^i \partial_j\dot{u}^i dx+\int_\Omega h\partial_i v^k \partial_k u^j\partial_j \dot{u}^idx\\&+\int_\Omega h\partial_j v^k  \left(\partial_i u^j+\partial_j u^i\right) \partial_k \dot{u}^idx\\
		\le &C\left\|\na\dot{u}\right\|_{L^2}\left\|\na v\right\|_{L^4}\left\|\na u\right\|_{L^4}\\
		\leq& \frac{\underline{\mu}}{4} \left\|\nabla \dot{u}\right\|_{L^2}^2+C\left\|\na{u}\right \|_{L^4}^4\\\leq&\frac{\underline{\mu}}{4} \left\|\nabla \dot{u}\right\|_{L^2}^2+C\left\|\sqrt{\rho}\dot{u}\right\|_{L^2}^2\left\|\na u\right\|_{L^2}^2.
	\end{aligned}
	\label{329}
\end{equation}

Moreover,  utilizing  periodic integration by parts and the fact $$\left(\partial_t+v\cdot\nabla \right)\partial_iP=\partial_i \dot P-\partial_i v\cdot\nabla P,$$ gives
\be 	 \ba
{J_4}=&-\int_\Omega \partial_{ i}\dot{P}\dot{u}^idx+ \int_\Omega   \partial_i v\cdot  \nabla P\dot{u}^i dx \\
=&\int_\Omega\left(P_t+v\cdot \nabla P\right)\partial_i v^k \partial_k u^i dx - \int_\Omega  P\partial_i v^k \partial_k \dot{u}^i   dx \\
=&\frac{d}{d t} \int_\Omega P \partial_{i} v^{k} \partial_{k} u^{i} dx-\int_\Omega P \left(\partial_{t}+v\cdot\nabla\right)\left( \partial_i v^k \partial_k {u}^i\right)dx - \int_\Omega P \partial_i v^k \partial_k \dot{u}^i   dx \\
=&\frac{d}{d t} \int_\Omega P \partial_{i} v^{k} \partial_{k} u^{i} dx  +\int_\Omega P \partial_{i} v^{j} \partial_{j} v^{k} \partial_{k} u^{i} dx+\int_\Omega P \partial_{i} v^{k} \partial_{k} v^{j} \partial_{j} u^{i} dx\\&-2 \int_\Omega P \partial_{i} v^{k} \partial_{k} \dot{u}^{i} dx-\int_\Omega P \partial_{k} {u}^{i}\partial_{i} \dot{v}^{k}  dx\\\triangleq&\frac{d}{d t} \int_\Omega P \partial_{i} v^{k} \partial_{k} u^{i}dx+\sum_{i=1}^{4}J_{4,i}.
\ea\ee
We estimate each $J_{4,i} (i=1,2,3,4)$ as follows.

First, since $\operatorname{div}v=0$, the $2\times2$ matrix $A = \na v$ has trace zero. We find $A^2=\lambda I$, where $\lambda=-\operatorname{det}(\na v)$. Combining this with $\operatorname{div}u=0$ gives
$$\partial_i v^j \partial_j v^k \partial_k u^i=\partial_i v^k \partial_k v^j \partial_j u^i=\operatorname{tr}\left(A^2\na u\right)=\lambda\operatorname{tr}(\na u)=\lambda\operatorname{div}u=0,$$
which implies that \begin{equation}
    \begin{aligned}\la{j41}
J_{4,1}+J_{4,2}=\int_\Omega P \partial_{i} v^{j} \partial_{j} v^{k} \partial_{k} u^{i} dx+\int_\Omega P \partial_{i} v^{k} \partial_{k} v^{j} \partial_{j} u^{i} dx=0.
    \end{aligned}
\end{equation}

Next,  by  $\operatorname{div}\left(\partial_{ i}v\right)=\partial_{ i}\operatorname{div}v, \nabla^{\perp} \cdot\left(\nabla \dot{u}^{i}\right) =0,$ one can deduce from  \eqref{ljo8}-\eqref{212}, \eqref{210} (for $p=4$),  \eqref{471}, and Lemma \ref{lemma2.3} that
\begin{equation}
    \begin{aligned}\la{j43}
\left|J_{4,3}\right|
        \leq&2\left|\int_\Omega P_1\partial_i v^k\partial_k\dot u^i\,dx\right|
+2\left|\int_\Omega P_2\partial_i v^k\partial_k\dot u^i\,dx\right|\\
\leq&C\|\nabla P_1\|_2\|\nabla v\|_2\|\nabla\dot u\|_2
+C\|P_2\|_4\|\nabla v\|_4\|\nabla\dot u\|_2\\
\leq&\frac{\underline\mu}{8}\|\nabla\dot u\|_2^2
+C\|\sqrt\rho\dot u\|_2^2\|\nabla u\|_2^2+C\|\nabla u\|_4^4\\
\leq&\frac{\underline\mu}{8}\|\nabla\dot u\|_2^2
+C\|\sqrt\rho\dot u\|_2^2\|\nabla u\|_2^2.
    \end{aligned}
\end{equation}

Finally, for $J_{4,4}$, we make the following split:
\begin{equation}\nonumber
	\begin{aligned}
	\left|J_{4,4}\right|
\leq&\bigg|\int_{\Om}P\partial_{k}u^i\partial_{ i}(\dot{u}^k *\omega_{\eta})dx\bigg|+\bigg|\int_{\Om}P\partial_{k}u^i\partial_{ i}\left[v\cdot(\na u^k*\omega_{\eta})-(v\cdot\na u^k)*\omega_{\eta}\right]dx\bigg|.
	\end{aligned}	
\end{equation}
Using  $\operatorname{div}\left(\partial_{ k}u\right)=\partial_{ k}\operatorname{div}u, \nabla^{\perp} \cdot\left[\nabla (\dot{u}^k *\omega_{\eta})\right] =0,$  \eqref{ljo8}-\eqref{212}, \eqref{210} (for $p=4$),  \eqref{471}, and Lemma \ref{lemma2.3} imply
\begin{equation}\nonumber
    \begin{aligned}
\bigg|\int_{\Om}P&\partial_{k}u^i\partial_{ i}(\dot{u}^k *\omega_{\eta})dx\bigg|\\\leq&C\left\|{P}_{1}\right\|_{\text{BMO}}\left\|\partial_{k} {u}^{i} \partial_{i} \left(\dot{u}^{k}*\omega_\eta\right)\right\|_{\mathcal{H}^{1}}+C\left\|P_{2}\right\|_{L^{4}}\left\|\nabla u\right\|_{L^{4}}\left\|\nabla \left(\dot{u}^{k}*\omega_\eta\right)\right\|_{L^{2}}\\\leq&C\left\|\na P_1\right\|_{L^2}\left\|\na u\right\|_{L^2}\left\|\na\dot{u}\right\|_{L^2}+C\left\|\na u\right\|_{L^4}^2\left\|\na\dot{u}\right\|_{L^2}\\\leq&\frac{\underline{\mu}}{8}\left\|\nabla \dot{u}\right\|_{L^{2}}^{2}+C\left\|\sqrt{\rho}\dot{u}\right\|_{L^2}^2\left\|\na u\right\|_{L^2}^2.
    \end{aligned}
\end{equation}
From     Lemmas \ref{lemma2.3} and \ref{lm26}, one has
\begin{equation}\nonumber
    \begin{aligned}
\bigg|\int_{\Om}P&\partial_{k}u^i\partial_{ i}\left[v\cdot(\na u^k*\omega_{\eta})-(v\cdot\na u^k)*\omega_{\eta}\right]dx\bigg|\\\leq&C\left\|{P}_{1}\right\|_{\text{BMO}}\left\|\partial_{k} {u}^{i} \partial_{ i}\left[v\cdot(\na u^k*\omega_{\eta})-(v\cdot\na u^k)*\omega_{\eta}\right]\right\|_{\mathcal{H}^{1}}\\&+C\left\|P_{2}\right\|_{L^{4}}\left\|\nabla u\right\|_{L^{4}}\left\|\na\left[v\cdot(\na u^k*\omega_{\eta})-(v\cdot\na u^k)*\omega_{\eta}\right]\right\|_{L^{2}}\\\leq&C\left\|\na P_1\right\|_{L^2}\left\|\na u\right\|_{L^2}\left\|\na \left[v\cdot(\na u^k*\omega_{\eta})-(v\cdot\na u^k)*\omega_{\eta}\right]\right\|_{L^2}\\&+C\left\|\nabla u\right\|_{L^{4}}^3\|\na v\|_{L^4}\\\leq&C\left\|\sqrt{\rho}\dot{u}\right\|_{L^2}\left\|\na u\right\|_{L^2}\left\|\na v\right\|_{L^4}\left\|\na u\right\|_{L^4}+C\left\|\na u\right\|_{L^4}^3\left\|\na v\right\|_{L^4}\\\leq&C\left\|\sqrt{\rho}\dot{u}\right\|_{L^2}^2\left\|\na u\right\|_{L^2}^2.
    \end{aligned}
\end{equation}
Then,
\begin{equation}
	\begin{aligned}\label{3312}
	\left|J_{4,4}\right|	\leq&\frac{\underline{\mu}}{8}\left\|\nabla \dot{u}\right\|_{L^{2}}^{2}+C\left\|\sqrt{\rho}\dot{u}\right\|_{L^2}^2\left\|\na u\right\|_{L^2}^2.
	\end{aligned}	
\end{equation}

Substituting  \eqref{329}-\eqref{3312} into \eqref{328},  we deduce from \eqref{210},  \eqref{2101}, and \eqref{332} that
\begin{equation}
	\Psi^{\prime}(t)+\frac{\underline{\mu}}{2}\left\|\nabla \dot{u} \right\|_{L^{2}}^2\leq C\left\|\sqrt{\rho}\dot{u}\right\|_{L^2}^2\left\|\na u\right\|_{L^2}^2,
	\label{333}
\end{equation}
where
\bnn
\Psi(t) \triangleq \frac{1}{2} \int_\Omega \rho\left|\dot{u}\right|^{2} dx-\int_\Omega P \partial_{i} v^{k} \partial_{k} u^{i} dx
\enn
satisfies
\begin{equation}\label{1000}
	\frac{1}{4} \int_\Omega \rho\left|\dot{u}\right|^{2} dx - C\left\|\nabla u\right\|_{L^2}^4 \leq \Psi(t) \leq \int_\Omega \rho\left|\dot{u}\right|^{2} dx + C\left\|\nabla u\right\|_{L^2}^4,
\end{equation}
due to   \eqref{310}.

Multiplying \eqref{333} by $\varsigma^2 e^{\alpha t}$, and applying \eqref{30}, \eqref{31}, \eqref{36},  \eqref{1000}, and Gronwall's inequality 	yields
\begin{equation}\label{336}
	\sup _{t \in[0, T]} \varsigma^2 e^{\alpha t}\left\|\sqrt{\rho} \dot{u}\right\|_{L^{2}}^{2}+\int_{0}^{T} \varsigma^2e^{\alpha t} \left\|\nabla \dot{u}\right\|_{L^{2}}^{2} d t \leq C.
\end{equation}
This and \eqref{33}, \eqref{210}, \eqref{2101}, \eqref{36}, imply  for $p\in[4, p^*)$,
\begin{equation}\label{3261}
	\sup _{t \in[0, T]} \varsigma^2 e^{\alpha t}\left(\left\|\na u\right\|_{L^p}^2+\left\|P\right\|_{L^p}^2\right)\leq C,
\end{equation}
which, together with  \eqref{36}, \eqref{3361},  \eqref{336}, and  \eqref{3261}, shows \eqref{326}. Thus,  the proof of Lemma \ref{lm32} is completed.\hfill $\Box$

\subsection{Dirichlet Problem}\la{3.3}

For a bounded Lipschitz domain $\Omega$, the closure defining $H^1_{0,\sigma}$ equals
$\{u\in H^1_0(\Omega)^2:\operatorname{div}u=0\}$ by
\cite[Lemma 2.3.3, Section 2.3.2]{Pokorny2022}.

 Let
\[
 \widetilde u(x)=\begin{cases}u(x),&x\in\Omega,\\0,&x\notin\Omega.\end{cases}
\]
For $u\in H^1_{0,\sigma}$, then
\[
 \widetilde u\in H^1(\mathbb R^2),\qquad
 \partial_i\widetilde u=\widetilde{\partial_i u},\qquad
 \operatorname{div}\widetilde u=0
 \quad\hbox{in }\mathcal D'(\mathbb R^2),
\]
and we set $v=S_\eta^\Omega u$  by the Stokes problem on  $\Om$ with $\displaystyle\int_\Omega P_\eta\,dx=0$:
\begin{equation}
\label{41}
 \begin{cases}
 -\Delta v+\nabla P_\eta=-\Delta(\widetilde u*\omega_\eta)&\text{in }\Omega,\\
 \operatorname{div}v=0&\text{in }\Omega,\\
 v=0&\text{on }\partial\Omega.\qquad

 \end{cases}
\end{equation}
By the standard Stokes theory for linear systems (see \cite[Theorem IV.6.1 (a),(b)]{G_1998}), system  \eqref{41} admits a  unique smooth solution up to the  boundary for fixed $\eta$. Moreover,  the standard Stokes estimates yield
\[\na v  \rightarrow \na u \text { in } L^{2}\left(\Omega\right),  \text { as } \eta \rightarrow 0,\]
and for any $1<p<\infty$,
\begin{equation}\label{43}
 \|\nabla v\|_{L^p(\Omega)}+\|P_\eta\|_{L^p(\Omega)}
 \le C({p,\Omega})\|\nabla u\|_{L^p(\Omega)},
\end{equation}
where $C$ is independent of $\eta$. Specially,  when $\Om=B_N$,  by transformation: $x=Ny$, $C$ is also independent of $N$.

For each smooth approximation at a fixed positive time, $u=v=0$ on $\partial\Omega$, and hence
\begin{equation}\label{eta-trace}
 \dot u=u_t+v\cdot\nabla u\in H^1_0(\Omega)^2,
 \quad
 \dot v=v_t+v\cdot\nabla v\in H^1_0(\Omega)^2.
\end{equation}
The divergence-free property does not generally hold for them. Indeed,
\begin{equation}\la{eta-trace0}
    \operatorname{div}\dot v=\partial_jv^k\partial_kv^j,
 \quad
 \int_\Omega\operatorname{div}\dot v\,dx=0,
 \quad
 \|\operatorname{div}\dot v\|_2\le C\|\nabla v\|_4^2.
\end{equation}

For  $F\in L^2(\Omega)^{2\times2}$ and $g\in L^2(\Omega)$ with
$\int_\Omega g\,dx=0$, we consider the problem
\[
-\Delta w+\nabla\pi=\operatorname{div}F,\quad
\operatorname{div}w=g,\quad w|_{\partial\Omega}=0.
\]
By \cite[Lemma 2.3.2, Section 2.3.2]{Pokorny2022}, there is a linear
right inverse of the divergence operator $\mathcal B_\Omega:L^2_0(\Omega)\to H^1_0(\Omega)^2$, with $\operatorname{div}\mathcal B_\Omega g=g$ and
\begin{equation}\label{weak-divergence-lifting}
\|\nabla\mathcal B_\Omega g\|_2\le C(\Omega)\|g\|_2,
\end{equation}
 where
$L^2_0(\Omega)=\{g\in L^2(\Omega):\int_\Omega g=0\}$. Moreover,
Lemma 2.3.3 of \cite{Pokorny2022} gives $w-\mathcal B_\Omega g\in H^1_{0,\sigma}$.
Testing with this field yields
$\|\nabla(w-\mathcal B_\Omega g)\|_2\le\|F\|_2+\|\nabla\mathcal B_\Omega g\|_2$,
and hence
\begin{equation}\label{nonzero-divergence-estimate}
\|\nabla w\|_2\le C(\Omega)(\|F\|_2+\|g\|_2).
\end{equation}
On $B_N$, define the scaled lifting by
\[
\mathcal B_{B_N}g(x)=N(\mathcal B_{B_1}[g(N\cdot)])(x/N),
\]
then
\[
\|\nabla\mathcal B_{B_N}g\|_2
 \le CN\|g(N\cdot)\|_{L^2(B_1)}=C\|g\|_{L^2(B_N)}.
\]
Thus, \eqref{nonzero-divergence-estimate} holds uniformly  in $N$ on $B_N$.

Substituting  $v=S_\eta^\Omega u$ from \eqref{41} into the approximate system \eqref{44}, and using the same arguments as in Lemma \ref{lm 30}, we arrive at the following exponential decay-in-time rate for $\|\sqrt{\rho}u\|_{L^2}$.
\begin{lemma}\label{lm 301}  Let $(\rho,u,h,P)$ be a smooth solution to \eqref{44}, \eqref{46}, \eqref{450}. Then there exists a positive constant $C$ depending only on  $\bar\rho$, $\underline{\mu}$,  $\left\|u_0\right\|_{L^2}$,  and $\Omega$ such that
	\begin{equation}
		\sup _{t \in[0, T]} e^{\alpha_1 t} \left\|\sqrt\rho u\right\|_{L^2}^{2} +\int_{0}^{T} e^{\alpha_1 t}\left\|\nabla u\right\|_{L^{2}}^{2} d t \leq C,
		\label{343}
	\end{equation}
	where 	$\alpha_1$ is some positive constant depending  only on  $\bar\rho $, $\underline{\mu}$, and $\Omega$.
\end{lemma}

Now we present some fundamental  results on $\na u$ and $ P$ from condition \eqref{0}.
\begin{lemma}\label{lm23}
	Let  $\left(\rho, u,h, P\right)$ be a smooth solution to  \eqref{44}, \eqref{46}, \eqref{450} satisfying \eqref{0}.    There exists   a positive constant  $\varepsilon$  depending only  on $\Om$  such that  for every $p\in[4,p^*)$,
	\begin{equation}\label{27}
		\left\|\na u\right\|_{L^p}\le C\left\|\na u\right\|_{L^2}^{2/p}\left\|\rho \dot{u}\right\|_{L^2}^{1-2/p},
	\end{equation}
	and
	\begin{equation}\la{271}
		\left\|P\right\|_{L^p}\leq	C\left\|{\rho}\dot{ u}\right\|_{L^2}+C\left\|\na u\right\|_{L^2},
	\end{equation}
	where    $C$  depends on $p$, $\bar{\mu}$, $\underline{\mu}$, $\Om$, and  $p^*>4$ is determined below by the relative oscillation of $h_0$, with $p^*=\infty$ if  this oscillation is zero.  
\end{lemma}
\noindent{\it Proof.}   Choose the constant $c_0$ in \eqref{280}.
It satisfies $\underline\mu\le c_0\le\bar\mu$, and the transport of $h$ preserves
$\|h(t)/c_0-1\|_{L^\infty}=\|h_0/c_0-1\|_{L^\infty}=\inf\limits_{c>0}\left\| {h_0}/{c}-1\right\|_{L^\infty}$.
Since $u_1=0$ on $\partial\Omega$,
$\int_\Omega\partial_j u_1^i\,dx=\int_{\partial\Omega}u_1^i n_j\,dS=0$.
Thus Lemmas \ref{lm21}, \ref{lm22}, and \eqref{ljo8}-\eqref{212} give
\begin{equation}\label{2811}
\begin{aligned}
\|\nabla u\|_{L^p}
&\le C(p,\Omega)\|\nabla u_1\|_{L^2}^{2/p}
       \|\nabla^2u_1\|_{L^2}^{1-2/p}+\|\nabla u_2\|_{L^p}\\
&\le C(p,\Omega){c_0}^{-(1-2/p)}
 \bigl(1+C_1(2,\Omega)\|h_0/c_0-1\|_{L^\infty}\bigr)^{2/p}
 \|\nabla u\|_{L^2}^{2/p}\|\rho\dot u\|_{L^2}^{1-2/p}\\
&\quad+C(p,\Omega)\|h_0/c_0-1\|_{L^\infty}\|\nabla u\|_{L^p}\\&\leq C(p,\Om, \underline{\mu}, \bar{\mu})\left\|\na u\right\|_{L^2}^{2/p}\left\|\rho\dot{u}\right\|_{L^2}^{1-2/p}+ {C_1}(p)\inf_{c>0}\left\|\frac{h_0}{c}-1\right\|_{L^{\infty}}\left\|\na u\right\|_{L^p}.
\end{aligned}
\end{equation}
Here, the constant  $C_1(p)$, arising from the fixed linear Stokes map, may be chosen locally bounded on finite $p-$intervals.
Take the continuous, strictly increasing controlling function
\begin{equation*}
\widetilde C_1(p)=1+p+
 \int_2^{p+1}\sup_{2\le q\le s+1}C_1(q)\,ds,\qquad p\ge2.
\end{equation*}
Then $\widetilde C_1(p)\ge C_1(p)$ and $\widetilde C_1(p)\to\infty$, as $p\to \infty$. Choose
\[
\varepsilon=\frac1{2\widetilde C_1(4)},\qquad
p^*=\widetilde C_1^{-1}\left(
\frac{1}{2\inf\limits_{c>0}\left\|\frac{h_0}{c}-1\right\|_{L^\infty}}\right).
\]
Specifically,  set $p^*=\infty$, when $\dis\inf_{c>0}\left\|\frac {h_0}{c}-1\right\|_{L^\infty}=0$.
Equivalently, we define
$\widetilde C(s)=\widetilde C_1^{-1}(1/(2s))$ for $0<s<\varepsilon$, with  $\widetilde C(0)=\infty$. Then $ \widetilde C$ is continuous and strictly decreasing on $(0,\varepsilon)$. So, for $4\le p<p^*$, $\dis {C_1}(p)\inf_{c>0}\left\|\frac{h_0}{c}-1\right\|_{L^{\infty}}<\frac12$ in \eqref{2811}   ,  hence
\begin{equation*}
\|\nabla u\|_{L^p}\le C\|\nabla u\|_{L^2}^{2/p}\|\rho\dot u\|_{L^2}^{1-2/p},
\end{equation*}
where  $C$ depends only on $p,\Omega,\underline\mu,\bar\mu$. On $B_N$, the same construction relies on the unit-ball Stokes constants and hence is independent of $N$. This,  combined with   \eqref{ljo8}-\eqref{212} yields that,  for $p\in[4, p^*)$,
\begin{equation*}\begin{aligned}
		\left\|P\right\|_{L^p}\leq& \left\|P_1\right\|_{L^p}+\left\|P_2\right\|_{L^p}\\
		\leq&C\left\|\na P_1\right\|_{L^2}+C\left\|\na u \right\|_{L^p}\\\leq
		&C\left\|\rho\dot{ u}\right\|_{L^2}+C\left\|\rho\dot{u}\right\|_{L^2}^{1-2/p}\left\|\na u\right\|_{L^2}^{2/p}\\\leq
		&C\left\|{\rho}\dot{ u}\right\|_{L^2}+C\left\|\na u\right\|_{L^2}.
	\end{aligned}
\end{equation*}
Thus, the proof of Lemma \ref{lm23} is finished.\hfill $\Box$

To estimate the pressure terms on  bounded domains, we  employ the following extension method,  since standard Hardy and BMO spaces do not apply  to  such domains.
\begin{lemma}\label{lemmaP1}
Let $f\in H^1_{0,\sigma}(\Omega)$, $g\in[H^1_0(\Omega)]^2$, and $P_1\in H^1(\Omega)$ with $\int_\Omega P_1\,dx=0$. Then
\begin{equation}\label{39a}
\left|\int_\Omega P_1\partial_i f^j\partial_j g^i\,dx\right|
\le C(\Omega)\|\nabla P_1\|_{L^2}\|\nabla f\|_{L^2}\|\nabla g\|_{L^2}.
\end{equation}
\end{lemma}
{\it Proof.}
We denote
$\tilde{f}$ and  $\tilde{g}$ as the zero extension of $f$ and $g$ to $\mr^2$, respectively, and $\tilde{P_1}$ as the general Sobolev extension  of $P_1$ (see \cite[Theorem 1, Section 5.4]{Evans_2010}),  satisfying
\begin{equation}
	\partial_i \tilde{f} =
	\left\{\begin{array}{l}
		\partial_i f, \qquad x\in \Omega,\\
		0, \qquad \quad x\in \mathbb{R}^2\setminus \Omega,
	\end{array}\right.\quad\partial_i \tilde{g} =
	\left\{\begin{array}{l}
		\partial_i g, \qquad x\in \Omega,\\
		0, \qquad \quad x\in \mathbb{R}^2\setminus \Omega,
	\end{array}\right.
	\label{348a}
\end{equation}
and
\begin{equation}
	\|\tilde{P}_{1}\|_{H^{1}\left(\mathbb{R}^{2}\right)} \leq C\left\|P_{1}\right\|_{H^{1}(\Omega)}.
\end{equation}
For the  ball $B_N$, let $P_1\in W^{1,q}(B_N)$ have mean zero, where $1<q<\infty$.
Applying the same extension theorem on $B_1$ to $P_1(N\cdot)$, yields an extension with support
in $B_2$, which is then scaled back to $\mathbb R^2$. Changing variables gives
\begin{equation}\label{52a}
\|\nabla\tilde P_1\|_{L^q(\mathbb R^2)}
\le C(q)\bigl(\|\nabla P_1\|_{L^q(B_N)}+N^{-1}\|P_1\|_{L^q(B_N)}\bigr)
\le C(q)\|\nabla P_1\|_{L^q(B_N)},
\end{equation}
due to  Poincar\'{e} inequality for a ball (see \cite[Theorem 2, Section 5.8.1]{Evans_2010}):
$$\left\|P_1\right\|_{L^q(B_N)}\leq C(q)N\left\|\na P_1\right\|_{L^q(B_N)}.$$

Then  recalling  $f|_{\partial\Omega}=0,$ it follows from  integration by parts that for any $\varphi \in \mathcal{D} (\mathbb{R}^2)$,
\begin{equation}
	-\int_{\mathbb{R}^{2}} \partial_{i} \tilde{f}^{j} \partial_{j} \varphi dx
	= \int_{\mathbb{R}^{2}} \tilde{f}^{j} \partial_{ij} \varphi dx= \int_{\Omega}  {f}^{j} \partial_{ij} \varphi dx=- \int_{\Omega} \partial_{j } {f}^{j} \partial_{ i} \varphi dx
	=0.
\end{equation}
This implies
\begin{equation} \label{315}	
	\mathrm{div} \partial_{j} \tilde{f} =0 \mbox{ in } \mathcal{D}'(\mathbb{R}^2),
\end{equation}
Similarly, one has
\begin{equation*}\nabla^{\perp} \cdot\left(\nabla \tilde{g}^{i}\right) =0 \mbox{ in } \mathcal{D}'(\mathbb{R}^2),\end{equation*}
which,  together with  \eqref{348a}-\eqref{315}, Lemma \ref{lemma2.3}, and  Poincar\'{e} inequality, gives that
\begin{equation*}
	\begin{aligned}	
		\left|\int_{\Omega} P_1 \partial_i f^j \partial_j g^i dx\right|& = \left|\int_{\mathbb{R}^2} \tilde{P_1} \partial_i { \tilde f}^j \partial_j \tilde{g}^i dx\right| \\
		&\leq C\|\tilde{P}_{1}\|_{\text{BMO}(\mathbb{R}^2)}\left\|\partial_{i} {f}^{j} \partial_{j} \tilde{g}^{i}\right\|_{\mathcal{H}^{1}(\mathbb{R}^2)} \\
		&\leq C\|\nabla \tilde{P}_{1}\|_{L^2(\mathbb{R}^2)}\|\nabla \tilde f\|_{L^{2}\left(\mathbb{R}^{2}\right)}\|\nabla \tilde{g}\|_{L^{2}\left(\mathbb{R}^{2}\right)}\\
		&\leq C\|\na P_{1}\|_{L^2(\Omega)} \|\nabla f\|_{L^{2}(\Omega)}\|\nabla g\|_{L^{2}(\Omega)}.
	\end{aligned}
\end{equation*}
This yields  \eqref{39a}. In view of \eqref{52a}, the
constant $C$ is independent of $N$ on $B_N$. Thus, the  proof of Lemma \ref{lemmaP1} is finished.
\hfill $\Box$

With Lemma \ref{lemmaP1} at hand, we establish   the exponential decay rate for $\|\nabla u\|_{L^2}.$
\begin{lemma}\label{lm3111}
	Let $(\rho,u,h,P)$ be a smooth solution to \eqref{44}, \eqref{46}, \eqref{450} satisfying \eqref{0}. For $p\in[4,p^*)$,  there exists  a positive constant $C$  depending only on $p$, $\bar{\mu}$, $\underline{\mu}$, $\bar\rho$, $\left\|u_0\right\|_{L^2}$,  and $\Omega$ such that
	\begin{equation}
		\sup _{t \in[0, T]}  \varsigma e^{\alpha_1 t}\left\|\nabla u\right\|_{L^{2}}^{2}+\int_{0}^{T}  \varsigma e^{\alpha_1 t}\left(\varsigma\left\|\sqrt{\rho}{u_t}\right\|_{L^2}^2+\left\|\sqrt{\rho}\dot{u}\right\|_{L^2}^2+\left\|\na u\right\|_{L^p}^2+\left\|P\right\|_{L^p}^2\right)dt \leq C,
		\label{345}
	\end{equation}
	where  $\alpha_1$ and $p^*$  are defined as  Lemmas \ref{lm 301} and \ref{lm23},  respectively.
\end{lemma}
\noindent{\it Proof.}
Multiplying $\eqref{44}_{3}$ by $\dot{u}$ and then integrating  over $\Omega$,  we find, after using  \eqref{eta-trace} and integration by parts
\begin{equation}
	\begin{aligned}
		\int_{\Omega} \rho|\dot{u}|^{2} dx=&\int_{\Omega}\left(\mathrm{div}\left(2hd\right)-\nabla P\right) \cdot \dot{u} dx\\
		=&\int_{\Omega}\mathrm{div}\left(2hd\right)\cdot \dot{u}dx+\int_{\Om}P\partial_{i}v^j\partial_{j}u^idx\\\triangleq&\sum_{i=1}^{2} K_{i} .
	\end{aligned}
	\label{346}
\end{equation}

We estimate $K_1, K_2$ as follows:

For $K_1$, the boundary term in the integration of \eqref{10000} vanishes  since $v=0$ on $\partial\Om$.
The estimate of $I_1$ implies, by \eqref{43} and \eqref{27} with $p=4$:
\begin{equation}\label{347}
K_1\le-\frac d{dt}\int_\Omega h|d|^2\,dx
+\frac14\|\sqrt\rho\dot u\|_2^2+C\|\nabla u\|_2^4.
\end{equation}

Then, for $K_2$,  by  \eqref{ljo8}-\eqref{212}, \eqref{27} (for $p=4$), \eqref{43},  and Lemma \ref{lemmaP1} with $f=v,g=u$, one checks that
\begin{equation}
	\begin{aligned}		
		|K_2|
&\leq\left|\int_\Omega P_1\partial_i v^j\partial_j u^i\,dx\right|
 +\left|\int_\Omega P_2\partial_i v^j\partial_j u^i\,dx\right|\\
&\le C\|\nabla P_1\|_{L^2}\|\nabla v\|_{L^2}\|\nabla u\|_{L^2}
 +\|P_2\|_{L^4}\|\nabla v\|_{L^2}\|\nabla u\|_{L^4}\\
&\le\frac14\|\sqrt\rho\dot u\|_{L^2}^2+C\|\nabla u\|_{L^2}^4.
	\end{aligned}
	\label{350}
\end{equation}

Substituting \eqref{347} and \eqref{350} into \eqref{346} gives
\begin{equation*}
	\frac{d}{d t} \int_{\Omega} h \left|d\right|^{2} dx+\frac12\int_{\Omega} \rho\left|\dot{u}\right|^{2} dx\le C\left\|\nabla u\right\|_{L^{2}}^4.
	\label{351}
\end{equation*}
Multiplying this by $\varsigma e^{\alpha_1 t}$, and by \eqref{343} and Gronwall's inequality, we  conclude that
\begin{equation*}
	\sup _{t \in[0, T]} \varsigma  e^{\alpha_1 t}\left\|\nabla u\right\|_{L^{2}}^{2}+\int_{0}^{T}  \varsigma e^{\alpha_1 t}\left\|\sqrt{\rho}\dot{u}\right\|_{L^2}^2\leq C.
\end{equation*}

Moreover, using  \eqref{471}, \eqref{31},  \eqref{41}$_3$ and Poincar\'{e} inequality  leads to
\begin{equation*}
	\begin{aligned}
		\int_{\Omega} \rho |u_{t}|^{2} d x & \leqslant \int_{\Omega} \rho|\dot{u}-v \cdot \nabla u|^{2} d x \\
		& \leq  C \left\|\sqrt{\rho} \dot{u}\right\|_{L^2}^2+C \int_{\Omega} \rho|v \cdot \nabla u|^{2} d x \\
		& \leq C \left\|\sqrt{\rho} \dot{u}\right\|_{L^2}^2+C\left \|v\right \|_{L^{4}}^2\left \|\nabla u\right \|_{L^{4}}^{2} \\
		& \leq C \left\|\sqrt{\rho} \dot{u}\right\|_{L^2}^2+C\left \|\nabla u\right \|_{L^{2}}^{2}\left\|\na u\right\|_{L^4}^2,
	\end{aligned}
\end{equation*}
which together with   \eqref{27} and \eqref{271} implies \eqref{345}. The proof of Lemma \ref{lm3111} is completed. \hfill $\Box$

Now,  we present the following essential connection between $\|\na \dot{v}\|_{L^2}$ and $\|\na\dot{u}\|_{L^2}$.

    \begin{lemma}\label{lm320}
Assume the assumptions of Lemma \ref{lm3111} hold. Then, on every positive-time interval,
\begin{equation}\label{53}
 \|\nabla\dot v\|_{L^2(\Omega)}
 \le C\|\nabla\dot u\|_{L^2(\Omega)}+C\|\nabla u\|_{L^4(\Omega)}^2,
\end{equation}
with  $C$ depending on $\Omega$ and  independent of $\eta$. On $B_N$, $C$ is also
independent of $N$. Moreover, for \eqref{41} on a fixed smooth domain, we have
\begin{equation}\label{weak-smoothing-L2-difference}
\|v-u\|_{L^2(\Omega)}\le C(\Omega)\sqrt\eta\|\nabla u\|_{L^2(\Omega)},
\qquad 0<\eta\le1.
\end{equation}
\end{lemma}
\noindent{\it Proof.} Applying $\partial_t+v\cdot\nabla$ to \eqref{41},  using \eqref{eta-trace0}, and
\[
 (\partial_t+v\cdot\nabla)\partial_i f
 =\partial_i((\partial_t+v\cdot\nabla)f)-\partial_i v^k\partial_k f,
\] gives
\begin{equation}\label{48}
\left\{
	\begin{aligned}
 -\Delta\dot v^j+\partial_j\dot P_\eta
=&-\Delta[(\partial_t+v\cdot\nabla)(\widetilde u^j*\omega_\eta)]+\partial_i[\partial_i v^k\partial_k(\widetilde u^j*\omega_\eta-v^j)]\\
 &+\partial_k[\partial_i v^k\partial_i(\widetilde u^j*\omega_\eta-v^j)]+\partial_k(P_\eta\partial_jv^k),\\
 \operatorname{div}\dot v=&\,\partial_jv^k\partial_kv^j,\\ \dot v|_{\partial\Omega}=&\,0.
\end{aligned}
	\right.
\end{equation}
Let $\widetilde v$ be the zero extension of $v$. Since $u_t,v$ have zero trace,
$\widetilde{\dot u}=\partial_t\widetilde u+\widetilde v\cdot\nabla\widetilde u$
in the sense of distributions.

Since \eqref{43} and Lemma \ref{lm26}, the first term on the right-hand side of \eqref{48} satisfies
\begin{equation}\label{49}
\begin{aligned}
\|\nabla[(\partial_t+v\cdot\nabla)(\widetilde u*\omega_\eta)]\|_{L^2(\Omega)}
&\le\|\nabla\widetilde{\dot u}*\omega_\eta\|_2
 +\|\nabla[\widetilde v\cdot\nabla(\widetilde u*\omega_\eta)
       -(\widetilde v\cdot\nabla\widetilde u)*\omega_\eta]\|_2\\
&\le\|\nabla\dot u\|_2+C\|\nabla v\|_4\|\nabla u\|_4
\le\|\nabla\dot u\|_2+C\|\nabla u\|_4^2,
\end{aligned}
\end{equation}
and by \eqref{43},
\begin{equation}\label{50}
 \begin{aligned}
 &\|\nabla v\|_4\|\nabla(\widetilde u*\omega_\eta-v)\|_4
 +\|P_\eta\|_4\|\nabla v\|_4
 +\|\nabla v\|_4^2\\
 &\qquad\le C\|\nabla u\|_4^2.
 \end{aligned}
\end{equation}
We apply \eqref{nonzero-divergence-estimate} to \eqref{48}, with  $F$ as the
velocity and pressure tensors  and
$g=\partial_jv^k\partial_kv^j$. Since  its mean is zero by \eqref{eta-trace0},
\eqref{49} and \eqref{50} imply \eqref{53}.

To prove \eqref{weak-smoothing-L2-difference}, for  $f\in L^2(\Omega)^2$ and let $(z, \tau) $ solve  the Stokes system
$$-\Delta z+\nabla\pi=f,\quad   \operatorname{div}z=0, \quad
z|_{\partial\Omega}=0,$$
with $\int_\Omega \pi=0$.  Then, by \cite[Theorem IV.6.1(a)]{G_1998},
$$\|z\|_{H^2}+\|\pi\|_{H^1}\le C(\Omega)\|f\|_2.$$
The weak form of \eqref{41} gives
\[
\int_\Omega(v-u)\cdot fdx
=\int_\Omega\nabla(\widetilde u*\omega_\eta-u):\nabla zdx
=\int_\Omega(\widetilde u*\omega_\eta-u)\cdot(-\Delta z)dx
 +\int_{\partial\Omega}(\widetilde u*\omega_\eta)\cdot\partial_nzdS.
\]
The first term is bounded by $C\eta\|\nabla u\|_2\|f\|_2$.
In a tubular neighborhood of the fixed smooth boundary, applying Jensen's inequality
and integrating along the normal direction, together with $u|_{\partial\Omega}=0$, yields
\[
\begin{aligned}
\|\widetilde u*\omega_\eta\|_{L^2(\partial\Omega)}^2
\le C(\Omega)\eta^{-1}
 \int_{\{x\in\Omega:\operatorname{dist}(x,\partial\Omega)<\eta\}}|u|^2\,dx\le C(\Omega)\eta\|\nabla u\|_2^2.
\end{aligned}
\]
Note that the estimate is valid for small $\eta$ and the remaining compact range
of $0<\eta\le1$ is absorbed into $C(\Omega)$ by the trace bound.

Since $\|\partial_nz\|_{L^2(\partial\Omega)}\le C(\Omega)\|f\|_2$,
$|\int_\Omega (v-u)\cdot f|\le C(\Omega)\sqrt\eta\|\nabla u\|_2\|f\|_2$.
Taking the supremum over $\|f\|_2\le1$ proves \eqref{weak-smoothing-L2-difference}.This completes the proof of  Lemma \ref{lm320}.\hfill $\Box$

With Lemma \ref{lm320} at hand, we  can  achieve  higher regularity for $\na u$ by following an approach similar to that in   Lemma \ref{lm32}.
\begin{lemma}\label{lm321}
	Let $(\rho,u,h, P)$ be a smooth solution to \eqref{44}, \eqref{46}, \eqref{450} satisfying \eqref{0}. For any $p\in[4,p^*)$,	  there exists a positive constant $C$ depending only on  $p$, $\bar\mu$, $\underline{\mu}$, $\bar\rho$, $\left\|u_0\right\|_{L^2}$,      and $\Omega$  such that
	\begin{equation}
		\sup _{t \in[0, T]} \varsigma^2 e^{\alpha_1t}\left(\varsigma^{\frac12}\left\|\sqrt{\rho} {u_t}\right\|_{L^{2}}^{2}+\left\|\sqrt{\rho} \dot{u}\right\|_{L^{2}}^{2}+\left\|\na u\right\|_{L^p}^2+\left\|P\right\|_{L^p}^2\right)+\int_{0}^{T} \varsigma^2 e^{\alpha_1t} \left\|\nabla \dot{u}\right\|_{L^{2}}^{2} d t \leq C,
		\label{352}
	\end{equation}
	where $\alpha_1$ and $p^*$ are the same as Lemmas \ref{lm 301} and  \ref{lm23}, respectively.
\end{lemma}
\noindent{\it Proof.} Operating $\left(\partial_t +v\cdot \nabla \right)$ to $\eqref{44}^{i}_3$ and integrating the resulting quality  multiplied by $\dot{u}^i$ over $\Omega$ give
\begin{equation}
	\begin{aligned}
		\frac{1}{2} &\frac{d}{dt} \int_\Omega \rho\left|\dot{u}\right|^{2} dx+ 2 \int_{\Omega} h\left|D\right|^{2} dx\\
		&=-\int_\Omega\dot{u}^i\partial_j\left(h\partial_j v^k \partial_k u^i\right)  dx-\int_\Omega \dot{u}^i\partial_j\left(h\partial_i v^k \partial_k u^j\right)  dx\\
		& \;\;\;\;-\int_\Omega \dot{u}^i\partial_j v^k \partial_k \left(h(\partial_i u^j+\partial_j u^i)\right) dx-\int_\Om\dot{u}^i \left(\partial_t +v\cdot \nabla \right)\partial_{ i}Pdx \\
		&\triangleq \sum_{i=1}^4 L_i.
	\end{aligned}
	\label{328111}
\end{equation}
In view of \eqref{43} and \eqref{27}, we have
\begin{equation}\label{L123-dirichlet}
\begin{aligned}
L_1+L_2+L_3&=\int_\Omega h\partial_jv^k\partial_ku^i\partial_j\dot u^i\,dx
 +\int_\Omega h\partial_iv^k\partial_ku^j\partial_j\dot u^i\,dx\\
&\quad+\int_\Omega h\partial_jv^k
       (\partial_i u^j+\partial_j u^i)\partial_k\dot u^i\,dx\\
&\le C\|\nabla\dot u\|_2\|\nabla v\|_4\|\nabla u\|_4\\
&\le\frac{\underline\mu}{4}\|\nabla\dot u\|_2^2
 +C\|\nabla u\|_2^2\|\sqrt\rho\dot u\|_2^2 .
\end{aligned}
\end{equation}
Proceeding as for  $J_4$, $J_{4,1}$, and $J_{4,2}$,  we arrive at
\begin{equation*}
\begin{aligned}
L_4
=&\frac d{dt}\int_\Omega P\partial_iv^k\partial_ku^i\,dx-2\int_\Omega P\partial_iv^k\partial_k\dot u^i\,dx
 \\
&-\int_\Omega P\partial_ku^i\partial_i\dot v^k\,dx .
\end{aligned}
\end{equation*}
Using   \eqref{ljo8}, \eqref{27}, \eqref{271},  \eqref{53}, and Lemma \ref{lemmaP1} with $f=v, g=\dot u$, yields
\begin{equation}\label{L4-mixed}
\begin{aligned}
2\left|\int_\Omega P\partial_iv^k\partial_k\dot u^i\,dx\right|
&\le C\left(\|\nabla P_1\|_2\|\nabla v\|_2
 +\|P_2\|_4\|\nabla v\|_4\right)\|\nabla\dot u\|_2\\
&\le\frac{\underline\mu}{4}\|\nabla\dot u\|_2^2
 +C\|\nabla u\|_2^2\|\sqrt\rho\dot u\|_2^2 .
\end{aligned}
\end{equation}
Similarly,
\begin{equation}\label{33121}
    \begin{aligned}
 \left|\int_\Omega P\partial_k u^i\partial_i\dot v^kdx\right|
 &\le C\left(\|\nabla P_1\|_2\|\nabla u\|_2
       +\|P_2\|_4\|\nabla u\|_4\right)\|\nabla\dot v\|_2\\
 &\le C\|\nabla u\|_2\|\sqrt\rho\dot u\|_2
 \left(\|\nabla\dot u\|_2+\|\nabla u\|_4^2\right)\\
&\le\frac{\underline\mu}{4}\|\nabla\dot u\|_2^2
 +C\|\nabla u\|_2^2\|\sqrt\rho\dot u\|_2^2.
\end{aligned}
\end{equation}

Putting \eqref{L123-dirichlet}-\eqref{33121}  into \eqref{328111}, implies
\begin{equation*}
 H'(t)+\frac{\underline\mu}{4}\|\nabla\dot u\|_2^2
 \le C\|\nabla u\|_2^2\|\sqrt\rho\dot u\|_2^2 ,
\end{equation*}where  \begin{equation*}
H(t)=\frac12\int_\Omega\rho|\dot u|^2\,dx
       -\int_\Omega P\partial_iv^k\partial_ku^i\,dx,
\end{equation*} satisfies
\begin{equation}
\frac14\|\sqrt\rho\dot u\|_2^2-C\|\nabla u\|_2^4
 \le H(t)\le
\|\sqrt\rho\dot u\|_2^2+C\|\nabla u\|_2^4,\label{100011}
\end{equation}
due to \eqref{350}.

Multiplying this by $\varsigma^2 e^{\alpha_1 t}$ and applying  \eqref{100011}, \eqref{345}, and  Gronwall's inequality, we  deduce that
\begin{equation*}
	\sup _{t \in[0, T]} \varsigma^2 e^{\alpha_1 t}\left\|\sqrt{\rho} \dot{u}\right\|_{L^{2}}^{2}+\int_{0}^{T} \varsigma^2 e^{\alpha_1 t} \left\|\nabla \dot{u}\right\|_{L^{2}}^{2} d t \leq C,
\end{equation*}
which, together with \eqref{3361}, \eqref{343}-\eqref{271},  gives  for $p\in[4, p^*)$,
\begin{equation*}
	\sup _{t \in[0, T]} \varsigma^2 e^{\alpha_1 t}\left(\varsigma^{\frac12}\left\|\sqrt\rho u_t\right\|_{L^2}^2+\left\|\na u\right\|_{L^p}^2+\left\|P\right\|_{L^p}^2\right)\leq C.
\end{equation*}
Hence, the proof of Lemma \ref{lm321} is completed.\hfill $\Box$

\subsection{Cauchy Problem}\label{3.2}

For the case of $\mr^2$, we establish a priori estimates uniform in $N$ under Dirichlet boundary condition \eqref{450} with $\Om={B_N}$. Here,  $v$ is defined as  described in \eqref{41}  with $\Om=B_N$.

First, without loss of generality, we assume
$$
\int_{\mr^2} \rho_0 d x=1.
$$
So, there  exists  a positive constant
$N_0$ such that
\begin{equation}\label{2711}
	\int_{B_{N_0}} \rho_0 d x \geq \frac{3}{4} \int_{\mr^2} \rho_0 d x=\frac{3}{4}.	
\end{equation}

The following key integrability properties hold for  $\na u$ and $P$.
\begin{lemma}\label{lm231}
	Let  $\left(\rho, u, h,P\right)$ be a smooth solution of  \eqref{44}, \eqref{46}, \eqref{450} satisfying \eqref{0}. Then  there exists  a positive constant  $\varepsilon$  depending  on the spatial dimensions   such that  for $p\in[4,p^*)$
\begin{equation}\label{272}
	\|\nabla u\|_{L^p(B_N)}
	\le C\|\nabla u\|_{L^2(B_N)}^{2/p}
	\|\rho\dot u\|_{L^2(B_N)}^{1-2/p},
\end{equation}
\begin{equation}\label{2712}
	\|\nabla u\|_{L^p(B_N)}+\|P\|_{L^p(B_N)}
	\le C\|\rho\|_{L^{p/2}(B_N)}^{1/2}
	\|\sqrt\rho\dot u\|_{L^2(B_N)} ,
\end{equation}
	where    $C$ depends only on    $p$, $\bar{\mu}$, $\underline{\mu}$,  and $p^*>4$ is determined below by the relative oscillation of $h_0$, with $p^*=\infty$ if  this oscillation is zero.
\end{lemma}
\noindent{\it Proof.} Take  the constant $c_0$ in \eqref{280}. Following the same argument as  in Lemma \ref{lm23} with $\Om$ replaced by $B_N$,  we  obtain the existence of  $\varepsilon$ in \eqref{0},  $p^*$ in \eqref{272} and
\begin{equation}\label{2720}
	\|\nabla u\|_{L^p(B_N)}
	\le C(p)\|\nabla u\|_{L^2(B_N)}^{2/p}\|\rho\dot u\|_{L^2(B_N)}^{1-2/p}
	+\bar{C}_1(p)\|h_0/{c_0}-1\|_{L^\infty}\|\nabla u\|_{L^p(B_N)} ,
\end{equation}
where the constants $\bar{C}_1(p)$ do not
depend on $N$. Choose $p^*>4$ so that
$\bar{C}_1(p)\|h_0/c_0-1\|_{L^\infty(B_N)}<1$ for each $4\le p<p^*$; when $h_0$ is
constant, take $p^*=\infty$. Then \eqref{2720} gives
\[
\bigl(1-\bar{C}_1(p)\|h_0/c_0-1\|_{L^\infty(B_N)}\bigr)\|\nabla u\|_{L^p(B_N)}
\le C(p)\|\nabla u\|_{L^2(B_N)}^{2/p}\|\rho\dot u\|_{L^2(B_N)}^{1-2/p},
\]
which proves \eqref{272}.

For the stronger estimate, every component of
$\nabla u_1$ has zero mean because $u_1|_{\partial B_N}=0$.
The Sobolev inequality and \eqref{211}, give
\begin{equation}\label{mass-stokes-sobolev}
	\begin{aligned}
		\|\nabla u_1\|_{L^p(B_N)}
		&\le C(p)\|\nabla^2u_1\|_{L^{2p/(p+2)}(B_N)}
		\le C(p)\|\rho\dot u\|_{L^{2p/(p+2)}(B_N)}\\
		&\le C(p)\|\sqrt\rho\|_{L^p(B_N)}\|\sqrt\rho\dot u\|_{L^2(B_N)}
		=C(p)\|\rho\|_{L^{p/2}(B_N)}^{1/2}\|\sqrt\rho\dot u\|_{L^2(B_N)} .
	\end{aligned}
\end{equation}
Using \eqref{212} for $u_2$ and absorbing as in \eqref{2720}, we obtain
\[
\bigl(1-\bar{C}_1(p)\|h_0/c_0-1\|_{L^\infty(B_N)}\bigr)\|\nabla u\|_{L^p(B_N)}
\le C_p\|\rho\|_{L^{p/2}(B_N)}^{1/2}\|\sqrt\rho\dot u\|_{L^2(B_N)} .
\]
Since $P_1$ has zero mean,  \eqref{ljo8}-\eqref{212} and \eqref{52a} give
\[
\begin{aligned}
	\|P\|_{L^p(B_N)}
	&\le C(p)\|\nabla P_1\|_{L^{2p/(p+2)}(B_N)}+C(p)\|(h-c_0)d(u)\|_{L^p(B_N)}\\
	&\le C(p)\|\rho\dot u\|_{L^{2p/(p+2)}(B_N)}
	+C(p)\|h-c_0\|_{L^\infty(B_N)}\|\nabla u\|_{L^p(B_N)}
	\\&\le C(p)\|\rho\|_{L^{p/2}(B_N)}^{1/2}\|\sqrt\rho\dot u\|_{L^2(B_N)},
\end{aligned}
\]
which, together with \eqref{2720}, completes the proof of Lemma \ref{lm231}.\hfill $\Box$

We now state the following key estimate for  $\left\|\na u\right\|_{L^2(B_N)}$, which is uniform in both time and
$N$.
\begin{lemma}\label{lm311}
	Let $(\rho,u,h,P)$ be a smooth solution to \eqref{44}, \eqref{46}, \eqref{450} satisfying \eqref{1101}. Then for $p\in[4,p^*)$
	\begin{equation}\label{337}
		\begin{aligned}
			\sup _{t \in[0, T]} &t\left\|\nabla u\right\|_{L^{2}(B_N)}^{2}+\int_{0}^{T}  t\left(\left\|\sqrt{\rho} \dot{u}\right\|_{L^{2}(B_N)}^{2}+ \left\|\na u\right\|_{L^p(B_N)}^2+ \left\|P\right\|_{L^p(B_N)}^2\right)dt \leq C,
		\end{aligned}
	\end{equation}
	where   $C$  depends only on $p$, $\bar{\mu}$, $\underline{\mu}$, $\left\|\rho_0\right\|_{L^1(B_N)\cap L^\infty(B_N)}$, and  $\left\|u_0\right\|_{L^2(B_N)}$.
\end{lemma}
\noindent{\it Proof.} First,  standard arguments (\cite[Theorem 2.1, Section 2]{Lions_1996}) show that
\begin{equation*}
	\left\|\rho(t)\right\|_{L^{1}\cap L^\infty(B_N)}=\left\|\rho_0\right\|_{L^{1}\cap L^{\infty}(B_N)},
\end{equation*}
and
\begin{equation}\label{339}
	\sup _{t \in[0, T]}\|\sqrt{\rho} u\|_{L^2(B_N)}^2+\int_0^T\|\nabla u\|_{L^2(B_N)}^2 d t \leq C.
\end{equation}

Then,  multiplying $\eqref{41}_3$ by $\dot{ u}$ and integrating  over $B_N$, we deduce from \eqref{346} with $\Om=B_N$ that
\begin{equation}\label{374}
	\begin{aligned}
		\int_{B_N} \rho|\dot{u}|^{2} dx
		=&\int_{B_N}\mathrm{div}\left(2hd\right)\cdot\dot{u}dx+\int_{B_N}P_2\partial_{i}v^j\partial_{j}u^idx+\int_{B_N}P_1\partial_{i}v^j\partial_{j}u^idx\\\triangleq&\sum_{i=1}^{3} O_{i} .
	\end{aligned}
\end{equation}
For $O_1, O_2$, applying the procedures  in \eqref{347} and \eqref{350} with $\Om=B_N$ leads to
\begin{equation}\label{375}
	O_1+O_2\leq -\frac{d}{d t} \int_{B_N} h|d|^{2} dx+\frac14\|\sqrt{\rho} \dot{u}\|_{L^{2}(B_N)}^{2}+C\|\nabla u\|_{L^{2}(B_N)}^{4}.
\end{equation}
Additionally, by \eqref{43}, \eqref{27}, \eqref{39a},   and Lemma \ref{lm22} with $\Om=B_N$,  we get
\begin{equation}\begin{aligned}
		|O_3|
&\le C\|\nabla P_1\|_{L^2(B_N)}\|\nabla v\|_{L^2(B_N)}
              \|\nabla u\|_{L^2(B_N)}\\
&\le\frac14\|\sqrt\rho\dot u\|_{L^2(B_N)}^2
 +C\|\nabla u\|_{L^2(B_N)}^4.
	\end{aligned}
	\label{39aa}
\end{equation}

Substituting \eqref{375} and \eqref{39aa}  into \eqref{374} gives
\begin{equation*}
	\frac{d}{d t} \int_{B_N} h \left|d\right|^{2} dx+\frac12\int_{B_N} \rho\left|\dot{u}\right|^{2} dx\le C\left\|\nabla u\right\|_{L^{2}(B_N)}^4.
\end{equation*}
Multiplying this by $t$, which  together with \eqref{339} and Gronwall's inequality, implies 
\begin{equation}\label{ball-low-core}
	\sup_{0<t\le T}t\|\nabla u(t)\|_{L^2(B_N)}^2
	+\int_0^Tt\|\sqrt\rho\dot u\|_{L^2(B_N)}^2dt\le C .
\end{equation}
By \eqref{2712},
\[
\begin{aligned}
	\int_0^Tt\bigl(\|\nabla u\|_{L^p(B_N)}^2+\|P\|_{L^p(B_N)}^2\bigr)\,dt
	&\le C\|\rho_0\|_{L^1}^{2/p}\|\rho_0\|_{L^\infty(B_N)}^{1-2/p}
	\int_0^Tt\|\sqrt\rho\dot u\|_{L^2(B_N)}^2\,dt
	\le C .
\end{aligned}
\]
Combining this with \eqref{ball-low-core} proves \eqref{337}. \hfill $\Box$

Following the approach as in Lemma \ref{lm321} with $\Om=B_N$, we establish  higher regularity for $\na u$.
\begin{lemma}\label{lm322}
	Let $(\rho,u,h,P)$ be a smooth solution to \eqref{44}, \eqref{46}, \eqref{450} satisfying \eqref{0}. 	Then   for any $p\in[4,p^*)$,
	\begin{equation}\begin{aligned}
			\sup _{t \in[0, T]} &t^2\left(\left\|\sqrt{\rho} \dot{u}\right\|_{L^{2}(B_N)}^{2}+ \left\|\na u\right\|_{L^p(B_N)}^2+ \left\|P\right\|_{L^p(B_N)}^2\right)\\&+\int_{0}^{T}t ^2 \left\|\nabla \dot{u}\right\|_{L^{2}(B_N)}^{2}  d t \leq C,
		\end{aligned}
		\label{340}
	\end{equation}
	where  $p^*$ is as  Lemma \ref{lm231} and   $C$ depends only on   $p$, $\bar\mu$, $\underline{\mu}$,  $\left\|\rho_0\right\|_{L^1(B_N)\cap L^\infty(B_N)}$, and $\left\|u_0\right\|_{L^2(B_N)}$.
\end{lemma}
\noindent{\it Proof.}
Similar to the proof of Lemma \ref{lm320}  and  \eqref{328111}-\eqref{100011}   with $\Om=B_N$, we have
\begin{equation}\label{391}
	\Theta^{\prime}(t)+\frac{\underline{\mu}}{4}\left\|\nabla \dot{u} \right\|_{L^{2}(B_N)}^2\leq C\left\|\sqrt{\rho}\dot{u}\right\|_{L^{2}(B_N)}^{2}\left\|\nabla u\right\|_{L^{2}(B_N)}^{2},
\end{equation}
where
\begin{equation*}
	\Theta \triangleq\frac{1}{2} \int_{B_N} \rho\left|\dot{u}\right|^{2} dx-\int_{B_N } P \partial_{i} v^{k} \partial_{k} u^{i}dx,
\end{equation*}
satisfies
\begin{equation}\label{1000111}
	\frac{1}{8} \int_{B_N}\rho|\dot{u}|^{2} dx - C\left\|\nabla u\right\|_{L^{2}(B_N)}^4 \leq \Theta(t) \leq \int_{B_N} \rho|\dot{u}|^{2} dx + C\left\|\nabla u\right\|_{L^{2}(B_N)}^4.
\end{equation}
Then, multiplying \eqref{391}  by $t^2$ and using \eqref{1000111} give
\begin{equation}\begin{aligned}\la{387}
		\frac{d}{dt}&\left(t^2\Theta(t)\right)+\frac{\underline{\mu}}{4}t^2\left\|\nabla \dot{u} \right\|_{L^{2}(B_N)}^2\\&\leq Ct^2\left\|\sqrt{\rho}\dot{u}\right\|_{L^{2}(B_N)}^{2}\left\|\na{u}\right\|_{L^2(B_N)}^2+Ct\left\|\sqrt{\rho}\dot{u}\right\|_{L^{2}(B_N)}^{2}+C\left\|\na u\right\|_{L^2(B_N)}^2.
	\end{aligned}	
\end{equation}
Applying \eqref{337}, \eqref{339}, \eqref{1000111}, \eqref{387}, and  Gronwall's inequality, we find
\begin{equation}\label{341}
	\sup _{t \in[0, T]} t^2 \left\|\sqrt{\rho} \dot{u}\right\|_{L^{2}(B_N)}^{2}+\int_{0}^{T} t^2  \left\|\nabla \dot{u}\right\|_{L^{2}(B_N)}^{2} d t \leq C,
\end{equation}
which, along with \eqref{272}, \eqref{2712}, and  \eqref{337} yields that for $p\in[4,p^*)$
\begin{equation}\label{342}
	\sup _{t \in[0, T]} t^{2}\left(\left\|\na u\right\|_{L^p(B_N)}^2+\left\|P\right\|_{L^p(B_N)}^2\right)\leq C.
\end{equation}
Thus, it follows from  \eqref{341} and   \eqref{342}  that   \eqref{340} holds.
This completes the  proof of Lemma \ref{lm322}.\hfill $\Box$

\section{Proofs of Theorems  \ref{tm1}-\ref{tm2}}\label{4}
With all the a priori estimates in Section \ref{3} at hand, we are now ready to prove Theorems \ref{tm1}-\ref{tm2}.

First, we have the following global existence result for  smooth solutions, which can be shown via Galerkin and transport method of   \cite[Theorem 2.4, Theorem 2.6]{Lions_1996} with  the standard theory of linear parabolic equations \cite{Evans_2010, LSU_1968}. 
\begin{pro}\label{4.1}
	For any given numbers $\eta$, $\bar{\rho}>0$, $\underline{\mu}>0$, and $\bar{\mu}>0$, assume that the smooth initial data $\left(\rho^\eta_0,   u^\eta_0, h^\eta_0, \right)$ satisfy
	$$
	0<\min\limits_{\overline{\Omega}}  \rho_0^\eta \leq \rho_0^\eta \leq 2\bar{\rho},\quad \underline{\mu} \leq h_0^\eta \leq \bar{\mu},\quad u_0^\eta\in C_{0, \sigma}^{\infty}\left(\Omega\right) .
	$$
	Then, the initial boundary value problem \eqref{44}, \eqref{46}, \eqref{45} or \eqref{44}, \eqref{46}, \eqref{450} admits a unique smooth solution $\left(\rho,  u, h, P\right)$ with $\int_\Om P\,dx=0$, satisfying
	$$
	\rho, u, h, P\in C^{\infty}\left(\overline{\Omega} \times(0, \infty)\right).
	$$
\end{pro}

\noindent{\bf Proof of Theorem \ref{tm1}.} The proof is divided into two steps.

    {\it Step 1.} We  construct approximated solutions to system \eqref{44}, \eqref{45} satisfying \eqref{periodic-conservation}  in the periodic  case.

For all $T>0$, and the initial data $(\rho_{0}, u_0)$ as described in \eqref{17},  we consider the approximated system \eqref{44}, \eqref{45} on $\Omega \times[0, T]$ with the initial conditions
\begin{equation}\label{4600}
	\rho(x, 0)=\rho^\eta_0(x), \quad  u(x, 0)= u^\eta_0(x), \quad h(x,0)= h^\eta_0(x),
\end{equation}
where  for $\eta>0,$ \begin{equation}\label{periodic-initial-normalization}
\begin{aligned}
&\rho_0^\eta=\rho_0*\omega_\eta+\eta,
\quad u_0^\eta=u_0*\omega_\eta
-\frac{\displaystyle\int_\Omega\rho_0^\eta(u_0*\omega_\eta)\,dx}
       {\displaystyle\int_\Omega\rho_0^\eta\,dx},\quad h_0^\eta=\mu(\rho_0)*\omega_\eta,
\end{aligned}
\end{equation}with \begin{equation}\label{periodic-initial-momentum}
\int_\Omega\rho_0^\eta u_0^\eta\,dx=0,
\quad\int_\Omega\rho_0^\eta\,dx
=\int_\Omega\rho_0\,dx+\eta|\Omega|
\ge\int_\Omega\rho_0\,dx>0.
\end{equation} Then,  for suitable $\eta>0$ and $\varepsilon_1$ described in \eqref{110},  $(\rho_0^\eta(x),{u}_0^\eta(x),h_0^\eta(x))$ satisfies
\begin{equation}\label{47}
	\begin{cases}
\left(\rho_0^\eta(x),u_0^\eta(x),h^\eta_0(x)\right)\in C^{\infty}_{\rm per}({\Om})\times C^\infty_{{\rm per},\sigma}(\Om)\times C^{\infty}_{\rm per}({\Om}),\\
		\eta\leq\rho^\eta_0\leq 2\bar{\rho}, \,\underline{\mu}\leq h^\eta_0\leq\bar{\mu},\,{\inf\limits_{c>0}\left\|\frac{h_0^\eta}{c}-1\right\|_{L^{\infty}(\Om)}}<\varepsilon_1,\\
		\left\|(\rho_{0}^\eta, h_{0}^\eta)\right\|_{L^\infty(\Om)}\leq2\left\|(\rho_{0}, \mu_0)\right\|_{L^\infty(\Om)}, \,\left\|  u_0^\eta\right\|_{L^2(\Om)}\leq C\left\|u_0\right\|_{L^2(\Om)},\,\operatorname{div} u^\eta_0=0,\\
		\lim\limits_{\substack{\eta \rightarrow 0}}\left(\left\|\rho_0^{\eta}-\rho_0\right\|_{L^1(\Om)}+\left\|u_0^{\eta}-\left(u_0-\frac{\int_\Omega\rho_0u_0\,dx}{\int_\Omega\rho_0\,dx}\right)\right\|_{L^2(\Om)}+\left\|h_0^{ \eta}-\mu_0\right\|_{L^1(\Om)}\right)=0.
	\end{cases}
\end{equation}
Then, by  Proposition \ref{4.1}, for all $T>0$, there exists a smooth solution $(\rho^\eta, u^\eta, h^\eta, P^\eta)$ to  \eqref{44}, \eqref{45}, \eqref{4600} on $\Om \times [0, T]$, with  $(\rho_{0}^\eta, u_{0}^\eta, h_0^\eta)$   satisfying \eqref{47} and $\int_\Om P^\eta dx=0$.  Integrating the mass and momentum equations in \eqref{44},  gives
\begin{equation}\label{periodic-approx-momentum}
\int_\Omega\rho^\eta(t)\,dx=\int_\Omega\rho_0^\eta\,dx,
\qquad\int_\Omega\rho^\eta u^\eta\,dx=0,
\qquad\int_\Omega\rho^\eta\dot u^\eta\,dx=0.
\end{equation}
Moreover,  \eqref{30},  \eqref{36}, and \eqref{326} hold with $(\rho, u, h, P)$ replaced by $(\rho^{\eta}, u^{\eta}, h^\eta, P^{\eta})$ for   some positive constant $C$ independent of $\eta$.

{\it Step 2.} We use the a priori estimates  from
Section \ref{3.1}  to establish the
existence and regularity of global weak solutions.

Indeed,  by passing to a subsequence (still denoted by $\eta$) and applying standard compactness arguments (see \cite[Section 2.3]{Lions_1996}),  \eqref{30}, \eqref{36}, \eqref{326}, and \eqref{47}, we conclude that the periodic problem \eqref{11}-\eqref{15} possesses a global weak solution $(\rho,  u, h, P)$ (with $\int_\Om P dx=0$)    on $\Om \times [0, T]$ satisfying for any $0<\tau<T<\infty$, as $\eta\rightarrow0$
\begin{equation}\la{412}
	\rho^{\eta} \rightharpoonup \rho \text { weakly-$*$} \text { in } L^{\infty}\left(\Om\times(0, T)\right),
\end{equation}
\begin{equation}\la{413}
	\rho^{\eta} \rightarrow \rho \text { in } C\left(\left[0, T\right] ; L^p\left(\Omega\right)\right), \forall 1\leq p <\infty,
\end{equation}
\begin{equation}\la{414}
	h^{\eta} \rightarrow h \text { in } C\left(\left[0, T\right] ; L^p\left(\Omega\right)\right), \forall 1\leq p <\infty,
\end{equation}\begin{equation}
	v^{\eta}-u^\eta \rightarrow 0 \text { strongly in } L^2\left(\Om\times(0, T)\right),
\end{equation}
\begin{equation}\la{44311110}
\hat P(\rho^\eta u^\eta)\longrightarrow\hat P(\rho u)
\quad\text{in }C([0,T];L^2_w(\Om)),
\end{equation}
\begin{equation}\label{compact-kinetic}
\sqrt{\rho^\eta}u^\eta\to\sqrt{\rho}u
\quad\text{strongly in }L^2\bigl(\Om\times(0,T)\bigr),
\end{equation}
\begin{equation}
    \sqrt{\rho^\eta}\dot u^\eta\rightharpoonup\sqrt\rho\dot u,
\text { weakly  in } L^2(\tau,T;L^2).
\end{equation}
\begin{equation}\la{4151}
	\nabla u^{\eta} \rightharpoonup \nabla u \text { weakly-$*$ in } L^{\infty}\left(\tau, T ; L^2\left(\Omega\right)\right),
\end{equation}
\begin{equation}\la{415}
	(\nabla u^{\eta},  P^{\eta} ) \rightharpoonup (\nabla u, P) \text { weakly in } L^{2}\left(\tau, T ; L^p\left(\Omega\right)\right), \forall   p\in[4, p^*),
\end{equation}
\begin{equation}\la{416}
(\nabla u^{\eta},  P^{\eta} )  \rightharpoonup  (\na u, P) \text { weakly-$*$ in } L^\infty\left(\tau,T; L^p(\Om)\right), \forall   p\in[4, p^*).
\end{equation}
Combining standard arguments (see \cite[Section 2.3]{Lions_1996}), with \eqref{412}-\eqref{416} yields $h(x, t)=\mu(\rho), \text{ a.e. in }\Om\times[0,T]$. This,  together with   a priori estimates in Section \ref{3.1}, proves the existence of global weak solution $(\rho,  u, P)$ (with $\int_\Om P dx=0$ and zero total momentum).  Then, by the Galilean invariance of the fluid-mechanics model under consideration, we can then recover  the original total momentum.  In particular, by \eqref{lj01} and \eqref{lj02}, we have
\begin{equation*}
\begin{aligned}
\nabla\big(P+R_iR_j(2h d_{ij})\big)=-Q(\rho\dot u),\quad
\partial_a(\hat P\otimes Q(2hd))_{ij}=\partial_a\partial_j\Delta^{-1}\hat P_{ik}(\rho\dot u)^k
 =-R_aR_j\hat P_{ik}(\rho\dot u)^k,
\end{aligned}
\end{equation*}
which implies \begin{equation}
\begin{aligned}\la{00}
\|\nabla(P+R_iR_j(2h d_{ij}))\|_2
 +\|\nabla(\hat P\otimes Q(2hd))\|_2\le C\|\rho\dot u\|_2
 \le C\sqrt{2\bar\rho}\,\|\sqrt\rho\dot u\|_2.
\end{aligned}
\end{equation}  Thus, the solution satisfies Definition \ref{d1} and the estimates \eqref{111}--\eqref{11110}. This completes the proof of Theorem \ref{tm1}.\hfill $\Box$

\noindent{\bf Proof of Theorem \ref{tm3}.} We make some modifications to the proof of Theorem \ref{tm1}.

For  the initial data $(\rho_{0}, u_0)$ as described in \eqref{17},   we set by \cite[Section 2.4]{Lions_1996},
\begin{equation}\la{533}
	\rho_{0}^\eta\triangleq\left(\rho_{0}\right)_\eta+\eta, \quad u^\eta_0\triangleq u_0^\eta,\quad h_0^\eta\triangleq\mu^\eta(\rho_0),
\end{equation}
and   for suitable $\eta$,   $(\rho_0^\eta,{u}_0^\eta,h_0^\eta)$ satisfies
\begin{equation}\label{470}
	\begin{cases}
		(\rho_0^\eta(x),u_0^\eta(x),h^\eta_0(x))\in C^{\infty}(\bar{\Om})\times C^\infty_{0,\sigma}(\Om)\times C^{\infty}(\bar{\Om}),\\
		\eta\leq\rho^\eta_0\leq 2\bar{\rho}, \,\,\underline{\mu}\leq h^\eta_0\leq\bar{\mu}, \,\, {\inf\limits_{c>0}\left\|\frac{h_0^\eta}{c}-1\right\|_{L^{\infty}(\Om)}}<\varepsilon_2,\\
		\left\|(\rho_{0}^\eta, h_{0}^\eta)\right\|_{L^\infty(\Om)}\leq2\left\|(\rho_{0}, \mu_0)\right\|_{L^\infty(\Om)}, \,\,\left\|u_0^\eta\right\|_{L^2(\Om)}\leq C \left\|u_0\right\|_{{L^2(\Om)}},\,\,\operatorname{div} u^\eta_0=0,\\
		\lim\limits_{\substack{\eta \rightarrow 0}}\left(\left\|\rho_0^{\eta}-\rho_0\right\|_{L^1(\Om)}+\left\|u_0^{\eta}-u_0\right\|_{L^2(\Om)}+\left\|h_0^{ \eta}-\mu_0\right\|_{L^1(\Om)}\right)=0.
	\end{cases}
\end{equation}
Now, for all $T>0$, we consider the  approximated system \eqref{44},  \eqref{450} on $\Omega \times[0, T]$ with the initial conditions \eqref{533} satisfying \eqref{470}.

By Proposition \ref{4.1}, there exists a smooth solution $(\rho^{\eta}, u^{\eta}, h^\eta, P^{\eta})$ to   \eqref{44},  \eqref{450}  on $\Om \times [0, T]$ with \eqref{533}. Furthermore, $(\rho^{\eta}, u^{\eta}, h^\eta, P^{\eta})$ satisfies Lemmas \ref{lm 301}-\ref{lm321}, where $C$ is independent of $\eta$. Thus, there exist some appropriate subsequences, still denoted by  $(\rho^{\eta}, u^{\eta}, h^\eta, P^{\eta})$,  such that  for any $0<\tau<T<\infty$, as  $\eta \rightarrow0$,
\begin{equation}\label{440}
	\rho^{\eta} \rightharpoonup \rho \text { weakly-$*$} \text { in } L^{\infty}\left(\Om\times(0, T)\right),
\end{equation}
\begin{equation}\la{441}
	\rho^{\eta} \rightarrow \rho \text { in } C\left(\left[0, T\right] ; L^p\left(\Omega\right)\right), \forall p\in[1,\infty),
\end{equation}
\begin{equation}\la{442}
	h^{\eta} \rightarrow h \text { in } C\left(\left[0, T\right] ; L^p\left(\Omega\right)\right), \forall   p\in[1,\infty),
\end{equation}
\begin{equation}
	v^{\eta}-u^\eta \rightarrow 0 \text { strongly in } L^2\left(\Om\times(0, T)\right),
\end{equation}
\begin{equation}\la{4431111}
\hat{P}({\rho^\eta} u^{\eta})  \rightarrow \hat{P}({\rho} u)  \text {  in } C\left([0, T] ; L^2_w\left(\Omega\right)\right),
\end{equation}
\begin{equation}\label{compact-kinetic1}
\sqrt{\rho^\eta}u^\eta\to\sqrt{\rho}u
\quad\text{strongly in }L^2\bigl(\Om\times(0,T)\bigr),
\end{equation}
\begin{equation}
    \sqrt{\rho^\eta}\dot u^\eta\rightharpoonup\sqrt\rho\dot u,
\text { weakly in } L^2(\tau,T;L^2).
\end{equation}
\begin{equation}\la{4431}
	\nabla u^{\eta} \rightharpoonup \nabla u \text { weakly-$*$ in } L^{\infty}\left(\tau, T ; L^2\left(\Omega\right)\right),
\end{equation}
\begin{equation}\la{443}
	(\nabla u^{\eta},  P^{\eta} )   \rightharpoonup (\nabla u, P) \text { weakly  in }L^{2}\left(\tau, T ; L^p\left(\Omega\right)\right), \forall p\in[4,p^*),
\end{equation}
\begin{equation}\la{445}
	(\nabla u^{\eta},  P^{\eta} )  \rightharpoonup (\na u, P) \text { weakly-$*$ in } L^\infty\left(\tau,T; L^p(\Om)\right), \forall p\in[4,p^*).
\end{equation}
By standard  arguments, together with \eqref{440}-\eqref{445} and  a priori estimates in Section \ref{3.3},  we deduce that $h(x, t)=\mu(\rho), \text{ a.e. in }\Om\times[0,T]$ and that $(\rho,  u, P)$ (with $\int_\Om P dx=0$)  is a global weak solution of \eqref{11}-\eqref{13}, \eqref{14}  satisfying \eqref{116}-\eqref{1161}. This completes the proof of Theorem \ref{tm3}. \hfill $\Box$

\noindent{\bf Proof of Theorem \ref{tm2}.}
For  $(\rho_{0}, u_0)$ as described in \eqref{17},  we consider the  approximated system \eqref{44},  \eqref{450} on $B_N \times[0, T]$ for any $T>0$ with initial conditions
\begin{equation}\la{5333}
	\rho_0^N \triangleq \rho_0*\omega_{N^{-1}}+N^{-1} e^{-|x|^2}, \quad u^N_0\triangleq {u}_{0}^N,\quad h_0^N \triangleq \mu\left(\rho_0\right)*\omega_{N^{-1}},
\end{equation}
where $N_0$ is given by
\eqref{2711}, and for suitably large $N$  and $\varepsilon_3$ as  in \eqref{1101},  $(\rho_0^N,{u}_0^N,h_0^N)$ satisfy
\begin{equation}\label{4701}
	\begin{cases}
		(\rho_0^N(x),u_0^N(x),h^N_0(x))\in C^{\infty}(\mr^2)\times C^\infty_{0,\sigma}(B_N)\times C^{\infty}(\mr^2),\\
		N^{-1} e^{-N^2}\leq\rho^N_0\leq	 2\bar{\rho}, \,\underline{\mu}\leq h^N_0\leq\bar{\mu},{\inf\limits_{c>0}\left\|\frac{h_0^N}{c}-1\right\|_{L^{\infty}(\mr^2)}}<\varepsilon_3,\, \int_{B_{N_0}} {\rho}_0^N d x\geq 1/2,\\
		\left\|(\rho_{0}^N, h_{0}^N)\right\|_{L^\infty(\mr^2)}\leq2\left\|(\rho_{0}, \mu_0)\right\|_{L^\infty(\mr^2)}, \left\|u_0^N\right\|_{L^2(\mr^2)}\leq C\left\|u_0\right\|_{L^2(\mr^2)},\,\operatorname{div} u^N_0=0,\\ \lim\limits_{\substack{N\rightarrow \infty}}\left(\left\|\rho_0^{N}-\rho_0\right\|_{L^1{(\mr^2)}}+\left\|u_0^{N}-u_0\right\|_{L^2(\mr^2)}+\left\|h_0^{N}-\mu_0\right\|_{L^1_{loc}(\mr^2)}\right)=0.
	\end{cases}
\end{equation}
Then,  for fixed $N>N_0$,  $0<\eta<1$, it follows from  Proposition \ref{4.1} that for any $T>0$, there exists a smooth solution $(\rho^{\eta,N}, u^{\eta,N}, h^{\eta,N}, P^{\eta,N})$ to   \eqref{44},  \eqref{450}, \eqref{5333} on $B_N \times [0, T]$.  Applying the same arguments  as  in Section \ref{3.2},  one easily checks that there exists some positive constant $C$, independent of $\eta$ and $N$,  such that   \eqref{337}, \eqref{339}, and \eqref{340} hold with $(\rho, u, h, P)$   replaced by $(\rho^{\eta,N}, u^{\eta,N},  h^{\eta,N}, P^{\eta,N})$.

Next, by passing to the limit first as $\eta\rightarrow 0$, then as $N\rightarrow\infty$ in  the approximate system \eqref{44}, \eqref{450}, \eqref{5333}, we claim that there exists a solution to  the initial problem \eqref{11}-\eqref{13}, \eqref{16}, \eqref{17},

Indeed, fixing  $N>N_0$,  passing to a subsequence as $\eta\rightarrow0$ (still denoted by $\eta$), we deduce that $(\rho^{\eta,N}, u^{\eta,N}, h^{\eta,N}, P^{\eta,N})$ converges weakly to some limit $(\rho^{N}, u^{N}, h^N,  P^{N})$, which satisfies
\eqref{337}, \eqref{339}, and \eqref{340}.
For every $0<\tau<T$, we have
\begin{equation}\label{ball-limit-traces}
u^N\in L^\infty(\tau,T;H^1_{0,\sigma}(B_N)),
\quad
\dot u^N=u_t^N+u^N\cdot\nabla u^N
       \in L^2(\tau,T;H^1_0(B_N)).
\end{equation}
Now, we  extend ($\rho^N, u^N, h^N$) by zero to $\mathbb{R}^2 \backslash B_N$ and denote the extension by $\left(\tilde{\rho}^N, \tilde{u}^N, \tilde h^N\right).$
Notice that
\begin{equation}
\label{ball-limit-dot-extension}
\begin{aligned}
\partial_t\widetilde u^N+
 \widetilde u^N\cdot\nabla\tilde u^N
 =\widetilde{\dot u^N}
 \quad\hbox{in }\mathcal D'(\mathbb R^2\times(\tau,T)),\quad \|\nabla\widetilde{\dot u^N}\|_{L^2(\mathbb R^2)}
       =\|\nabla\dot u^N\|_{L^2(B_N)} .
\end{aligned}
\end{equation}
 It follows from \eqref{52a}, \eqref{337}, \eqref{339}, \eqref{340},   and \eqref{4701}  that for all $T>0$, and any $p\in [4,p^*)$,
\begin{equation}\la{531a}
	\begin{aligned}
		\sup _{0 \leq t \leq T}&\left(\left\|\sqrt{\tilde{\rho}^{N}} \tilde{u}^{N}\right\|_{L^2}^2+\varsigma\left\|\nabla \tilde{u}^N\right\|_{L^2}^2+\varsigma^2\left\|\sqrt{\tilde\rho^N}\dot {\tilde{u}}^N\right\|_{L^2}^2+\varsigma^2\left\|\na \tilde u^N\right\|_{L^p}^2\right)\\&\leq \sup_{0 \leq t \leq T}\left(\left\|\sqrt{{\rho}^{N}} {u}^{N}\right\|_{L^2(B_N)}^2+\varsigma\left\|\nabla {u}^N\right\|_{L^2(B_N)}^2+\varsigma^2\left\|\sqrt{\rho^N}\dot u^N\right\|_{L^2(B_N)}^2 +\varsigma^2\left\|\na u^N\right\|_{L^p(B_N)}^2\right),
	\end{aligned}
\end{equation}
and \begin{equation}\label{54}
	\begin{aligned}
		\int_0^T&\left(\left\|\na \tilde u^N\right\|_{L^2}^2+\varsigma\left\|\sqrt{{\tilde\rho}^N}\dot{\tilde{u}}^N\right\|_{L^2}^2+\varsigma^2\left\|\na \dot{ \tilde{u}}^N\right\|_{L^2}^2+\varsigma\left\|\na \tilde u^N\right\|_{L^p}^2\right)dt\\&\leq \int_0^T\left(\left\|\na u^N\right\|_{L^2(B_N)}^2+\varsigma\left\|\sqrt{\rho^N}\dot{u}^N\right\|_{L^2(B_N)}^2+\varsigma^2\left\|\na \dot{ u}^N\right\|_{L^2(B_N)}^2+\varsigma\left\|\na u^N\right\|_{L^p(B_N)}^2\right)dt.
	\end{aligned}
\end{equation}
For fixed $T$, choose $\widetilde N>2N_0$ and a smooth cutoff
$\varphi_{\widetilde N}$ satisfying
\begin{equation}\label{ball-mass-cutoff}
0\le\varphi_{\widetilde N}\le1,\quad
\varphi_{\widetilde N}=1\text{ on }B_{\widetilde N/2},\quad
\operatorname{supp}\varphi_{\widetilde N}\subset B_{\widetilde N},\quad
|\nabla\varphi_{\widetilde N}|\le C/\widetilde N.
\end{equation}
Then, $\eqref{11}_1$ gives
\begin{equation}\label{ball-mass-cutoff-derivative}
\begin{aligned}
\left|\frac d{dt}\int_{B_N}\rho^N\varphi_{\widetilde N}\,dx\right|
&=\left|\int_{B_N}\rho^Nu^N\cdot\nabla\varphi_{\widetilde N}\,dx\right|\\
&\le \frac C{\widetilde N}
 \left(\int_{B_N}\rho^N\,dx\right)^{1/2}
 \|\sqrt{\rho^N}u^N\|_2
 \le \frac C{\widetilde N}.
\end{aligned}
\end{equation}
Choose $\widetilde N$  sufficiently large such that $CT/\widetilde N\le\frac14\int_{\mathbb R^2}\rho_0\,dx$. By \eqref{4701}, integration of \eqref{ball-mass-cutoff-derivative} gives
\begin{equation}\label{ball-mass-cutoff-lower}
\inf_{0\le t\le T}\int_{B_{\widetilde N}}\rho^N(t,x)\,dx
\ge\frac14\int_{\mathbb R^2}\rho_0\,dx.
\end{equation}

With \eqref{531a}, \eqref{54}, \eqref{ball-mass-cutoff-lower} at hand, we  obtain that there exist some appropriate subsequences, still denoted by  $(\tilde{\rho}^{N}, \tilde{u}^{N}, \tilde{h}^{N},  \tilde{P}^{N})$,  such that  for any $0<\tau<T<\infty$, as  $N\rightarrow\infty$,
\begin{equation}\label{4401}
	\tilde\rho^{N} \rightharpoonup \rho \text { weakly-$*$} \text { in } L^{\infty}\left(\Om\times(0, T)\right),
\end{equation}
\begin{equation}\la{4411}
	\tilde\rho^{N} \rightarrow \rho \text { in } C\left(\left[0, T\right] ; L^p\left(\Omega\right)\right), \forall   p\in[1, \infty),
\end{equation}
\begin{equation}\la{4421}
	\tilde h^{N} \rightarrow h \text { in } C\left(\left[0, T\right] ; L^p_{\text{loc }}\left(\Omega\right)\right), \forall   p\in[1, \infty),
\end{equation}
\begin{equation}\la{443111100}
\hat{P}({\tilde\rho^N} \tilde u^{N})  \rightarrow \hat{P}({\rho} u) \text {  in } C\left([0, T] ; L^2_w\left(\Omega\right)\right),
\end{equation}
\begin{equation}\la{h00}
\int_0^T\int_\Om\left(\tilde\rho^{N} \tilde u^{N} \otimes \tilde u^{N}(t)-\rho u \otimes u(t)\right) :\nabla \psi dx d t\rightarrow0, \quad \forall\psi\in\Psi,
\end{equation}
\begin{equation}\la{44311}
	\nabla \tilde u^{N} \rightharpoonup \nabla u \text { weakly-$*$ in } L^{\infty}\left(\tau, T ; L^2\left(\Omega\right)\right),
\end{equation}\begin{equation}
    \sqrt{\tilde\rho^\eta}\widetilde{\dot u^\eta}\rightharpoonup\sqrt\rho\dot u,
\text { weakly-$*$ in } L^2(\tau,T;L^2).
\end{equation}
\begin{equation}\la{44310000}
	\nabla \tilde u^{N} \rightharpoonup \nabla u \text { weakly-$*$ in } L^{2}\left(\tau, T ; L^p\left(\Omega\right)\right)\cap L^{\infty}\left(\tau, T ; L^p\left(\Omega\right)\right), \forall   p\in[4, p^*),
\end{equation}
where \eqref{h00} is proved in \cite[Section 4]{wei_2025}.
Finally, by standard  arguments along with \eqref{00} in the whole-space, \eqref{4401}-\eqref{44310000}, the   a priori estimates in Section \ref{3.2}, and Lemma 2.3 of \cite{Lu_2017}), we conclude that   $h(x, t)=\mu(\rho), \text{ a.e. in }\Om\times[0,T]$ and the limit $(\rho,  u, \na P)$ is a global weak solution of \eqref{11}-\eqref{13}, \eqref{16}  satisfying \eqref{1111}-\eqref{11611}.  Thus, the proof of Theorem \ref{tm2} is completed.  \hfill $\Box$

\end{document}